\documentclass{article}
\usepackage{amsmath,amssymb}
\usepackage[all]{xy}
\newtheorem{thm}{Theorem}[section]
\newtheorem{prop}[thm]{Proposition}
\newtheorem{cor}[thm]{Corollary}
\newtheorem{lem}[thm]{Lemma}
\newenvironment{dfn}{\medskip\refstepcounter{thm}
\noindent{\bf Definition \thesection.\arabic{thm}\ }}{\medskip}

\newenvironment{proof}[1][,]{\medskip\ifcat,#1
\noindent{{\it Proof}:\ }\else\noindent{\it Proof of #1.\ }\fi}
{\hfill$\square$\medskip}
\newenvironment{note}[1][Note]{\begin{trivlist}
\item[\hskip \labelsep {\bfseries #1}]}{\end{trivlist}}
\newenvironment{notes}[1][Notes]{\begin{trivlist}
\item[\hskip \labelsep {\bfseries #1}]}{\end{trivlist}}
\newenvironment{remark}[1][Remark]{\begin{trivlist}
\item[\hskip \labelsep {\bfseries #1}]}{\end{trivlist}}
\newenvironment{remarks}[1][Remarks]{\begin{trivlist}
\item[\hskip \labelsep {\bfseries #1}]}{\end{trivlist}}
\def\eq#1{{\rm(\ref{#1})}}
\def\H{{\mathbb H}}
\def\P{{\mathbb P}}

\def\R{{\mathbb R}}
\def\O{{\mathbb O}}

\def\N{{\mathbb N}}
\DeclareMathOperator\GL{GL}

\def\SU{\mathop{\rm SU}}

\def\U{\mathbin{\rm U}}

\def\d{{\rm d}}
\def\w{\wedge}
\def\C{{\mathbb C}}

\def\Re{\mathop{\rm Re}\nolimits}
\def\Im{\mathop{\rm Im}\nolimits}

\DeclareMathOperator{\Ker}{Ker}
\DeclareMathOperator{\Coker}{Coker}
\DeclareMathOperator{\cs}{cs}
\DeclareMathOperator{\dR}{dR}
\DeclareMathOperator{\GG}{G}
\DeclareMathOperator{\Aut}{Aut}
\DeclareMathOperator{\loc}{loc}
 
\begin{document}

\title{Deformation Theory of Asymptotically Conical Coassociative 4-folds}
\author{\textsc{Jason D. Lotay}\\ University College\\ Oxford}
\date{}
\maketitle \linespread{1.1}

\begin{center}
{\large\bf Abstract}
\end{center}

\noindent Suppose that a coassociative 4-fold $N$ in $\R^7$ is asymptotically conical to a cone $C$ with 
 rate $\lambda<1$.  If $\lambda\in[-2,1)$ is generic, we show that the moduli space of coassociative deformations of $N$ which 
are also asymptotically conical to $C$ with rate $\lambda$ is a smooth manifold, and we calculate its dimension.  If $\lambda<-2$ and generic, we show that the moduli space is locally homeomorphic to the kernel of a smooth map between 
smooth manifolds, and we give a lower bound for its expected dimension.  We also derive a test for when $N$ will be planar if 
$\lambda<-2$ and we discuss examples of asymptotically conical coassociative 4-folds.

\section{Introduction}

\footnote{\textsl{2000 Mathematics Subject Classification: 53C38.}}In this article we study \textit{coassociative} 4-folds $N$ in $\R^7$
which are \emph{asymptotically conical} (AC). 
 Our main result (cf.~Theorem \ref{modulispacethm} \& Corollary \ref{modulispacecor}) is the following.

\begin{thm}\label{mainthm}
Let $N$ be a coassociative 4-fold in $\R^7$ which is AC with rate $\lambda<1$ to a cone $C$.  Let 
$\mathcal{M}(N,\lambda)$ be the moduli space of coassociative deformations of $N$ which are also AC with 
rate $\lambda$ to $C$.
\begin{itemize}
\item[\emph{(a)}] For generic $\lambda\in[-2,1)$ the 
deformation theory is unobstructed so $\mathcal{M}(N,\lambda)$ is a smooth manifold near $N$.
\item[\emph{(b)}] For generic $\lambda<-2$, $\mathcal{M}(N,\lambda)$ is locally homeomorphic to the kernel of a 
smooth map between smooth manifolds.
\end{itemize}
\end{thm}

\noindent In case (a), the dimension of $\mathcal{M}(N,\lambda)$ is equal to the dimension of the \emph{infinitesimal deformation 
space}.  We determine this dimension explicitly in Proposition \ref{dimprop1} if additionally $\lambda\in[-2,0)$, and we give
 upper and lower bounds for it otherwise in Proposition \ref{dimprop2}.  
In case (b), the map can be considered as a projection from the infinitesimal deformation space to the
\emph{obstruction space} and hence $\mathcal{M}(N,\lambda)$ is smooth if the latter space is zero.  Moreover, we give a lower bound 
for the expected dimension of $\mathcal{M}(N,\lambda)$ for $\lambda<-2$ in Proposition \ref{dimprop3}.
These dimensions are given in terms of topological data from $N$ and $C$, and an analytic quantity determined by the cone $C$.


\smallskip

This article is motivated by the work of McLean \cite[$\S$4]{McLean}
on compact coassociative 4-folds and Marshall
\cite{Marshall} on deformations of AC special Lagrangian
(SL) submanifolds. The latter was in fact studied earlier by
Pacini \cite{Pacini} using different methods. Asymptotically
conical SL submanifolds are also discussed in the 
 papers
by Joyce \cite{Joyce5}-\cite{Joyce9} on SL submanifolds with 
conical singularities.  Deformations of asymptotically \emph{cylindrical} coassociative
4-folds are studied by Joyce and Salur in \cite{JoySalur} and Salur in \cite{Salur}.  Other coassociative deformation theories have been studied by the author: coassociative 4-folds with conical singularities in \cite{Lotaysing} and 
 compact coassociative 4-folds with boundary in \cite{LotayKov} (with Kovalev). 

\smallskip

We begin, in $\S$\ref{dfnsection}, by defining 
the submanifolds that we study here.  Much of our work 
is analytic in nature, so 
 to obtain many of
the results we use \emph{weighted} Sobolev spaces,
which are a natural choice when studying AC
submanifolds in this way.  Thus, in $\S$\ref{weightedsection}, we define the
weighted Banach spaces that we require. 

In $\S$\ref{defmapsection} we construct a \emph{deformation map} which corresponds locally to the moduli space
$\mathcal{M}(N,\lambda)$.  We may therefore view the kernel of the linearisation of the deformation map at zero as the 
infinitesimal deformation space.  We also define an associated map which
is \emph{elliptic} at zero since its derivative there acts as
$d+d^*$ from pairs of self-dual 2-forms and 4-forms to 3-forms.  This allows to prove that the kernel of the deformation 
map consists of smooth forms.

In $\S$\ref{dd*section} we discuss the Fredolm and index theory of the operator $d+d^*$: in particular, we describe a
 countable discrete
set $\mathcal{D}$, depending only on 
 $C$, 
consisting of the rates $\lambda$ for which $d+d^*$ is not Fredholm.   Our main result of the section (Theorem \ref{splitthm}) 
identifies the obstruction space for our deformation theory. 

In $\S$\ref{modulispacesection} we prove the deformation theory results which lead directly to Theorem \ref{mainthm}.  Section 
$\S$\ref{dimsection} then contains our aforementioned dimension calculations.

In $\S$\ref{invsection}, we construct two invariants of $N$ and hence derive a test for when $N$
will be planar.  Finally, in $\S$\ref{exsection}, we discuss examples of AC coassociative 4-folds, including explicit examples 
which have rate $-3/2$.  
  There are no known concrete examples of AC coassociative 4-folds with rate $\lambda<-2$, but such submanifolds are essential for 
  the desingularization theory of coassociative 4-folds with conical singularities, as discussed in \cite{Lotaydesing}.

\medskip

\noindent\textbf{Acknowledgements} \hspace{2pt} Many thanks are due
to Dominic Joyce for his hard work and guidance during the course of this
project.  I would also like to thank EPSRC for providing the funding for this study.

\section{Asymptotically conical coassociative\\ submanifolds of {\boldmath $\R^7$}}
\label{dfnsection}
 
We begin by defining \emph{coassociative 4-folds} in $\R^7$, for which we
introduce a distinguished 3-form 
on $\R^7$, following the notation
in \cite[Definition 11.1.1]{Joyce1}.

\begin{dfn}\label{phidfn} Let $(x_1,\ldots,x_7)$ be coordinates on
$\R^7$ and write
$d{\bf x}_{ij\ldots k}$ for the form $dx_i\w dx_j\w\ldots\w dx_k$.
Define a 3-form $\varphi$ by:
\begin{equation*}
\varphi = d{\bf x}_{123}+d{\bf x}_{145}+d{\bf
x}_{167}+d{\bf x}_{246}- d{\bf x}_{257}-d{\bf x}_{347}-d{\bf
x}_{356}.
\end{equation*}

\noindent The 4-form $\ast\varphi$, where $\varphi$ and
$\ast\varphi$ are related by the Hodge star, is given by:
\begin{equation*}
\ast\varphi = d{\bf x}_{4567}+d{\bf x}_{2367}+d{\bf
x}_{2345}+d{\bf x}_{1357}-d{\bf x}_{1346}-d{\bf x}_{1256}-d{\bf
x}_{1247}.
\end{equation*}
The subgroup of $\GL(7,\R)$ preserving the 3-form $\varphi$ is
$\GG_2$.
\end{dfn}

We can now characterise coassociative
4-folds in $\R^7$.

\begin{dfn}\label{coassdfn} A 4-dimensional submanifold $N$ of $\R^7$ is coassociative if and only if
$\varphi |_N\equiv 0$ and $*\varphi|_N>0$.
\end{dfn}

\noindent This is not the standard definition, which is formulated in terms of \emph{calibrated geometry}, 
but is equivalent to it by \cite[Proposition IV.4.5 \& Theorem IV.4.6]{HarLaw}.  

\begin{remark} The vanishing of $\varphi$ on a 4-fold $N$ in $\R^7$
forces $*\varphi$ to be nowhere vanishing on $N$.  Thus, the condition $*\varphi|_N>0$ amounts to a choice of orientation.
\end{remark}

In order that we may study deformations of coassociative 4-folds, we need an important elementary
result \cite[Proposition 4.2]{McLean}.

\begin{prop}\label{jmathprop}
Let $N$ be a coassociative 4-fold in $\R^7$. There is an
isomorphism $\jmath_N$ between the normal bundle $\nu(N)$ of $N$ in $\R^7$ and
$\Lambda^2_+T^*N$ given by $\jmath_N:v\mapsto (v\cdot\varphi)|_{TN}$.
\end{prop}

\begin{note}  
Let $u\in T\R^7$.  There exist unique normal and tangent vectors $v$ and $w$ on $N$ such that
$u|_N=v+w$.  Since $\varphi$ vanishes on the coassociative 4-fold $N$, $(w\cdot\varphi)|_{TN}=0$.  Thus,
$\jmath_N(u)=\jmath_N(v)\in\Lambda^2_+T^*N$.  So, $\jmath_N:T\R^7|_N\rightarrow\Lambda^2_+T^*N$.
\end{note}

Before 
defining AC submanifolds we 
clarify what we mean by a cone in $\R^n$.

\begin{dfn}\label{conedfn}
A \emph{cone} $C\subseteq\R^n$ is a nonsingular submanifold, except
perhaps at 0, satisfying $e^tC=C$ for all $t\in\R$.  We call $\Sigma=C\cap\mathcal{S}^{n-1}$
the \emph{link} of $C$.  
\end{dfn}

\begin{dfn}
\label{ACsubmflddfn} Let $C$ be a closed cone in $\R^n$, let $\Sigma$ be the link of $C$ and let $N$
be a closed submanifold of $\R^n$.  Then $N$ is
\emph{asymptotically conical} (AC) \emph{to $C$} (with rate $\lambda$) if there exist
constants $\lambda<1$ and $R>1$, a
compact subset $K$ of $N$, and a diffeomorphism
$\Psi:(R,\infty)\times\Sigma\rightarrow N\setminus K$ such that
\begin{equation}
\label{Psieq1} \big|\nabla^j\big(\Psi(r,\sigma)-\iota(r,\sigma)\big)\big|=O\big(r^{\lambda-j}\big) 
\qquad\text{for $j\in\N$ as $r\rightarrow\infty$},
\end{equation}
 \noindent where $\iota:(0,\infty)\times\Sigma\rightarrow \R^n$ is 
 the inclusion map 
given by $\iota(r,\sigma)=r\sigma$. Here $|\,.\,|$ is calculated using the conical metric
$g_{\text{cone}}=dr^2+r^2g_{\Sigma}$ on $(0,\infty)\times\Sigma$, where $g_{\Sigma}$ is the round
metric on $\mathcal{S}^{n-1}$ restricted to $\Sigma$, and $\nabla$ is a combination of the
Levi--Civita connection derived from $g_{\text{cone}}$ and the flat
connection on $\R^n$, which acts as partial differentiation.
\end{dfn}

We also make the following definition.

\begin{dfn}
\label{ACmflddfn} Let $(M,g)$ be a Riemannian $n$-manifold. 
Then $M$ is \emph{asymptotically conical} (AC) (with rate $\lambda$)
if there exist constants $\lambda<1$ and $R>1$, a compact $(n\!-\!1)$-dimensional Riemannian
submanifold $(\Sigma,g_{\Sigma})$ of $\mathcal{S}^{n-1}$, compact $K\subseteq M$, and a
diffeomorphism $\Psi:(R,\infty)\times\Sigma\rightarrow M\setminus
K$ such that
\begin{equation}
\label{Psieq2}
\big|\nabla^j\big(\Psi^*(g)-g_{\,\text{cone}}\big)\big|=O\big(r^{\lambda-1-j}\big)\qquad\text{for $j\in\N$ as $r\rightarrow\infty$,}
\end{equation}
where $(r,\sigma)$ are coordinates on $(0,\infty)\times\Sigma$,
$g_{\,\text{cone}}=dr^2+r^2g_{\Sigma}$ 
 on $(0,\infty)\times\Sigma$, 
$\nabla$ is the Levi--Civita connection of
 $g_{\,\text{cone}}$ and $|\,.\,|$ is calculated using
$g_{\text{cone}}$.

Define $\iota:(0,\infty)\times\Sigma\rightarrow\R^n$ by $\iota(r,\sigma)=r\sigma$ and let
$C=\Im\iota$
.  Then $C$ is the \emph{asymptotic cone} of $M$ and  
the components of $M_{\infty}=M\setminus K$ are the \emph{ends} of $M$.
\end{dfn}

\vspace{-8pt}

\noindent The condition $\lambda<1$ in the definition above ensures that, by \eq{Psieq2}, the metric $g$ on $M$ 
converges to $g_{\text{cone}}$ at infinity.  

\begin{note} By comparing \eq{Psieq1} and \eq{Psieq2}, we see
that if $N$ is a Riemannian submanifold of $\R^n$ which is AC to a cone $C$ with rate $\lambda$, it can be considered as an AC manifold
with rate $\lambda$ and asymptotic cone $C$. 
\end{note}

\begin{dfn}\label{radiusfndfn} Let $M$ be an AC manifold and use the notation of
Definition \ref{ACmflddfn}. A \emph{radius function}
$\rho:M\rightarrow [1,\infty)$ on $M$ is a smooth function such that
there exist positive constants $c_1<1$ and $c_2>1$ with
$c_1r<\Psi^*(\rho)<c_2r$. 
\end{dfn}

\noindent If $M$ is AC we may define a radius function $\rho$ 
by setting $\rho =1$ on $K$,
$\rho\big(\Psi(r,\sigma)\big)=r$ for $r>R+1$, and then
extending $\rho$ smoothly to our required function on $M$.

\medskip

We conclude the section with the following elementary result.
\begin{prop}\label{coassconeprop}
Suppose that $N$ is a coassociative 4-fold in $\R^7$ which is AC with rate $\lambda$ 
to a cone $C$ in $\R^7$.  Then
$C$ is coassociative.
\end{prop}
\begin{proof}
By Definition \ref{coassdfn}, we have that $\varphi|_N\equiv 0$. We
also have, using \eq{Psieq1}, that
$$\big|\Psi^*\big(\varphi|_N\big)-\varphi|_{C}\big|=O\big(r^{\lambda-1}\big)\qquad \text{as $
r\rightarrow\infty$.}$$ 
Therefore, since $\lambda-1<0$, $\big|\,\varphi|_{C}\big|\rightarrow 0$ as $r\rightarrow\infty$.  However,
$\big|\,\varphi|_{C}\big|$ must be independent of $r$ since
$T_{r\sigma}C=T_{\sigma}C$ for all $r>0$,
$\sigma\in\Sigma$. Hence $\varphi|_{C}\equiv 0$.
\end{proof}

\begin{notes}\begin{itemize}\item[]
\item[(a)] Manifolds are taken to be nonsingular and submanifolds to
be embedded, for convenience, unless stated otherwise.\vspace{-4pt}
\item[(b)] We use the convention that the natural numbers
$\N=\{0,1,2,\ldots\}$.
\end{itemize} \end{notes}

\section{Weighted Banach spaces}
\label{weightedsection}   

We shall define \emph{weighted} Banach spaces of forms on an AC manifold, following
\cite[$\S$1]{Bartnik}.  We use the notation and definition of the usual `unweighted' Banach spaces 
as  
in \cite[$\S$1.2]{Joyce1}; 
that is, Sobolev and H\"older spaces are denoted by $L^p_k$ and $C^{k,\,a}$ respectively,
where $p\geq 1$, $k\in\N$ and $a\in(0,1)$.  
We also introduce the notation $C^k_{\loc}$ for the space of forms $\xi$ such that $f\xi$ lies in $C^k$ for every
smooth compactly supported function $f$, and similarly define spaces $L^p_{k,\,\loc}$ and $C^{k,\,a}_{\loc}$.

For the whole of this section, we let $(M,g)$ be an AC $n$-manifold and $\rho$ be a radius function on $M$ as in
Definition \ref{radiusfndfn}.

\begin{dfn}\label{wSobdfn} 
%
%
Let $p\geq1$, $k\in\N$ and $\mu\in\R$.  
The \emph{weighted Sobolev space}
$L_{k,\,\mu}^p(\Lambda^mT^*M)$ of $m$-forms on $M$ is the 
subspace of $L_{k,\,\loc}^p(\Lambda^mT^*M)$ such that 
\begin{equation*}
\|\xi\|_{L_{k,\,\mu}^p}=\left(\sum_{j=0}^k\int_M
|\rho^{j-\mu}\nabla^j\xi|^p\rho^{-n} \,dV_g\right)^\frac{1}{p}
\end{equation*}
is finite.
The normed vector space $L_{k,\,\mu}^p(\Lambda^mT^*M)$ is a Banach space. 
\end{dfn}

\begin{note} 
\begin{equation}\label{Lpeq}
L^p(\Lambda^mT^*M)=L_{0,\,-\frac{n}{p}}^p(\Lambda^mT^*M).
\end{equation}
\end{note}

We now define \emph{dual} weighted Sobolev space which shall be invaluable later.

\begin{dfn}\label{dualdfn}  Use the notation from Definition \ref{wSobdfn}. 
Let $p,q>1$ such that $\frac{1}{p}+\frac{1}{q}=1$, let $k,l\in\N$ and let $\mu\in\R$.  Define
a pairing $\langle\,.\,,\,.\,\rangle:L^p_{k,\,\mu}(\Lambda^mT^*M)\times L^q_{l,\,-n-\mu}(\Lambda^mT^*M)\rightarrow\R$ by
$$\langle\xi,\eta\rangle=\|\xi\w\ast\eta\|_{L^1}.$$
We shall occasionally refer to this as the \emph{dual pairing}.  For our purposes, we take the dual space of
 $L^p_{k,\,\mu}(\Lambda^mT^*M)$ to be $L^q_{l,\,-n-\mu}(\Lambda^mT^*M)$, with linear functionals represented by dual pairings.
\end{dfn}


\begin{dfn}\label{wCkdfn} 
Let $\mu\in\R$ and let $k\in\N$.  The
\emph{weighted $C^k$-space} $C_{\mu}^{k}(\Lambda^mT^*M)$ of
$m$-forms on $M$ is the subspace of $C^k_{\loc}(\Lambda^mT^*M)$ 
such that the norm 
$$\|\xi\|_{C_{\mu}^{k}}=\sum_{j=0}^k \sup_M|\rho^{j-\mu}\nabla^j\xi|$$
is finite.  We also define
$C_{\mu}^{\infty}(\Lambda^mT^*M)=\cap_{k\geq
0}C_{\mu}^{k}(\Lambda^mT^*M)$.  Then $C_{\mu}^{k}(\Lambda^mT^*M)$
is a Banach space but in general $C_{\mu}^{\infty}(\Lambda^mT^*M)$
is not.  Notice that we have a continuous embedding $C^k_{\mu}\hookrightarrow C^l_{\nu}$ whenever $k\geq l$ and $\mu\leq\nu$.
\end{dfn}

We conclude this spate of definitions by defining weighted
\emph{H\"older} spaces.  These shall not be required for the
majority of the paper as weighted Sobolev and $C^k$-spaces will often suffice.

\begin{dfn}\label{wHolderdfn}
Let $a\in(0,1)$, $k\in\N$ and $\mu\in\R$.
Let $d(x,y)$ be
the geodesic distance between points $x,y\in M$, let $0<c_1<1<c_2$ be constant  
and let 
\begin{align*}
H=\{(x,y)\in M&\times M\,:\,x\neq
y,\,c_1\rho(x)\leq\rho(y)\leq c_2\rho(x)\,\;\text{and}\;\,\\
&\text{there exists a geodesic in $M$ of length $d(x,y)$ from $x$ to
$y$}\},
\end{align*} 
A section $s$
of a vector bundle $V$ on $M$, endowed with Euclidean metrics on its fibres and a connection preserving these metrics, is \emph{H\"older continuous} (with
\emph{exponent $a$}) if
$$[s]^a=
\sup_{(x,y)\in H}\frac{|s(x)-s(y)|_{V}}{d(x,y)^a}<\infty.$$ 
We understand the quantity $|s(x)-s(y)|_V$ as follows.  Given
$(x,y)\in H$, there exists a 
geodesic $\gamma$ of length $d(x,y)$ connecting $x$ and $y$.
Parallel translation along $\gamma$ using the connection on $V$
identifies the fibres over $x$ and $y$ and the metrics on them.
Thus, with this identification, $|s(x)-s(y)|_V$ is well-defined.

%
%
The \emph{weighted H\"older space}
$C_{\mu}^{k,\,a}(\Lambda^mT^*M)$ of $m$-forms $\xi$ on $M$ is the
subspace of $C^{k,\,a}_{\loc}(\Lambda^mT^*M)$ such
that the norm 
$$\|\xi\|_{C^{k,\,a}_{\mu}}=\|\xi\|_{C^{k}_{\mu}}+[\xi]^{k,\,a}_{\mu}$$
is finite, where
$$[\xi]^{k,\,a}_{\mu}=[\rho^{k+a-\mu}\nabla^k\xi]^{a}.
$$  Then $C_{\mu}^{k,\,a}(\Lambda^mT^*M)$ is a
Banach space.  It is clear that we have an embedding
$C_{\mu}^{k,\,a}(\Lambda^mT^*M)\hookrightarrow
C_{\mu}^l(\Lambda^mT^*M)$ whenever $l\leq k$.
\end{dfn}

\noindent The set $H$ in the definition above is introduced so that $[\xi]^{k,\,a}_{\mu}$ is well-defined.


We shall need the analogue of the Sobolev Embedding Theorem for
weighted spaces, which is adapted from \cite[Lemma
7.2]{LockhartMcOwen} and \cite[Theorem 1.2]{Bartnik}.
\begin{thm}[Weighted Sobolev Embedding Theorem]
\label{wSobthm} 
Let $p$, $q>
1$, $a\in (0,1)$, $\mu,\nu\in\R$ and $k,l\in\N$.
\begin{itemize}
\item[{\rm (a)}] If\/ $k\geq l$, $k-\frac{n}{p}\geq l-\frac{n}{q}$
and either
$p\leq q$ and $\mu\leq\nu$, or 
$p>q$ and $\mu<\nu$,
there is a continuous embedding
$L_{k,\,\mu}^p(\Lambda^mT^*M)\hookrightarrow$
$L_{l,\,\nu}^q(\Lambda^mT^*M)$. 
\item[{\rm (b)}] If\/
$k-\frac{n}{p}\geq l+a$, there is a continuous embedding
$L_{k,\,\mu}^p(\Lambda^mT^*M)\hookrightarrow
C_{\mu}^{l,\,a}(\Lambda^mT^*M)$.
\end{itemize}
\end{thm}

We shall also require an Implicit Function Theorem for Banach
spaces, which follows immediately from \cite[Theorem 2.1]{Lang2}.
\begin{thm}[Implicit Function Theorem]
\label{IFthm} $\!$Let $X$ and $Y\!$ be Banach spaces and let $U\subseteq
X$ be an open neighbourhood of $0$.  Let $\mathcal{F}:U\rightarrow Y$ be a
$C^k$-map $(k\geq 1)$ such that $\mathcal{F}(0)=0$.  Suppose further that
$d\mathcal{F}|_0:X\rightarrow Y$ is surjective with kernel $K$ such
that $X=K\oplus A$ for some closed subspace $A$ of $X$.  There
exist open sets $V\subseteq K$ and $W\subseteq A$, both containing
$0$, with $V\times W\subseteq U$, and a unique $C^k$-map
$\mathcal{G}:V\rightarrow W$ such that
$$\mathcal{F}^{-1}(0)\cap(V\times W)=\big\{\big(x,\mathcal{G}(x)\big)\,:\,x\in V\big\}$$
in $X=K\oplus A$.
\end{thm}

\section{The deformation map}
\label{defmapsection}
 
For the rest of the paper, let $N\subseteq\R^7$ be a coassociative 4-fold which is
asymptotically conical to a cone $C\subseteq\R^7$ with rate $\lambda$, 
and use the notation of Definition \ref{ACsubmflddfn}.   
In particular, $(0,\infty)\times\Sigma\,{\buildrel \iota\over\cong}\, C$, where
$\Sigma=C\cap\mathcal{S}^6$, with coordinates $(r,\sigma)$ on $(0,\infty)\times\Sigma$, and $(R,\infty)\times\Sigma
\,{\buildrel \Psi\over\cong}\, N\setminus K$, where $K$ is a compact subset of $N$ and $R>1$.  Moreover, let
$\rho$ be a radius function on $N$, as defined in Definition \ref{radiusfndfn}, and choose
$\Psi$
uniquely by imposing the condition that
\begin{equation}\label{Psiorthogeq}
\Psi(r,\sigma)-\iota(r,\sigma)\in (T_{r\sigma}C)^\perp\quad\text{for all
$(r,\sigma)\in(R,\infty)\times\Sigma$,}
\end{equation}
which can be achieved by
making $R$ and $K$ larger if necessary.

We shall endeavour to remind the reader of this notation when required.

\subsection{Preliminaries}

We wish to discuss deformations of $N$; that is, coassociative submanifolds of $\R^7$ which are `near' to $N$.  
We define this formally.

\begin{dfn}\label{modulispacedfn}
The \emph{moduli space of deformations} $\mathcal{M}(N,\lambda)$ is
the set of coassociative 4-folds $N^\prime\subseteq \R^7$ which are
AC to $C$ with rate $\lambda$ such that there exists a
diffeomorphism $h:\R^7\rightarrow \R^7$, with $h(N)=N^{\prime}$, isotopic to the identity.
\end{dfn}

The first result we need is immediate from the proof of
\cite[Chapter IV, Theorem 9]{Lang}.
\begin{thm}\label{tubenbdthm} Let $P$ be a closed submanifold
of a Riemannian manifold $M$.  There exist an open subset $V$ of the
normal bundle $\nu(P)$ of $P$ in $M$, containing the zero section,
and an open set $S$ in $M$ containing $P$, such that the exponential
map $\exp|_V:V\rightarrow S$ is a diffeomorphism.
\end{thm}

\begin{note} The proof of this result relies entirely on the
observation that $\exp|_{\nu(P)}$ is a local isomorphism upon the
zero section.
\end{note} 

\noindent This information helps us prove a useful corollary.

\begin{cor}\label{tubenbdcor}  Let
$P=\iota\big((R,\infty)\times\Sigma\big)$, $Q=N\setminus K$ and define
$n_P: \nu(P)\rightarrow\R^7$ by $n_P(r\sigma,v)=
v+\Psi(r,\sigma)$. 
There exist an open subset $V$ of $\nu(P)$, containing the
zero section, and an open set $S$ in $\R^7$ containing $Q$, such
that $n_P|_{V}:V\rightarrow S$ is a diffeomorphism. 
Moreover, $V$ can be chosen to be an open neighbourhood of the zero section in $C^1_1$.
\end{cor}

\begin{proof}
Note that $n_P$ takes the zero section of $\nu(P)$ to $Q$. By the
definition of $\Psi$ in Definition \ref{ACsubmflddfn}, $n_P$ is a local isomorphism upon the zero
section. Thus, by the note after Theorem \ref{tubenbdthm}, we have open sets
$V$ and $S$ such that $n_P|_V:V\rightarrow S$ is a diffeomorphism.

Since $\Psi-\iota$ is orthogonal to $(R,\infty)\times\Sigma$ by \eq{Psiorthogeq}, it can
be identified with a small section of the normal bundle. Hence $P$
lies in $S$ as long as $S$ grows with order $O(r)$ as $r\rightarrow
\infty$.  As we can form $S$ and $V$ in a translation equivariant
way because we are working on a portion of the cone $C$, we can
construct our sets with this growth rate at infinity and such that
they do not collapse near $R$.  Thus, we can ensure that $V$ is an open set in $C^1_1$.
\end{proof}

Recall 
that, by Proposition \ref{coassconeprop}, $C$ is coassociative.  
Therefore, since $\Psi(r,\sigma)-\iota(r,\sigma)$ lies in
$(T_{r\sigma}C)^{\perp}\cong\nu_{r\sigma}(C)$ for $r>R$ by \eq{Psiorthogeq},
$\Psi-\iota$ can be identified with the graph of an element
$\gamma_C$ of $\Lambda_+^2T^*P$ by
Proposition \ref{jmathprop}, using the notation of Corollary \ref{tubenbdcor}. Then
\begin{equation}\label{gammaCeq}
|\nabla_C^j\gamma_C|=O\big(r^{\lambda-j}\big)\quad\text{for}\;
j\in\N\;\text{as}\; r\rightarrow\infty, \end{equation} since $N$ is
AC to $C$ with rate $\lambda$, where $\nabla_C$ and $|\,.\,|$ are
the Levi--Civita connection and modulus calculated using the conical metric.  Thus, $\gamma_C$ lies in $
C_{\lambda}^{\infty}(\Lambda_+^2T^*P)$. Moreover, we have a
decomposition:
$$\R^7=T_{\Psi(r,\sigma)}N\oplus\nu_{r\sigma}(C)$$ at
$\Psi(r,\sigma)$.  We can therefore identify
$\nu_{\Psi(r,\sigma)}(N)$ with $\nu_{r\sigma}(C)$ and hence 
identify $\Lambda_+^2T^*N$ and $\Lambda_+^2T^*C$ near infinity.
Formally, we have the following.

\begin{prop}\label{Upsilonprop} Use the notation of Corollary
\ref{tubenbdcor} and let $\jmath_C$ and $\jmath_N$ be the
isomorphisms given by Proposition \ref{jmathprop} applied to
$C$ and $N$ respectively.  There exists a diffeomorphism
$\Upsilon:\nu(P)\rightarrow\nu(Q)$, with $\Upsilon(0)=0$, and hence
a diffeomorphism
$\tilde{\Upsilon}:\Lambda^2_+T^*P\rightarrow\Lambda^2_+T^*Q$ given
by $\tilde{\Upsilon}=\jmath_N\circ\Upsilon\circ\jmath_C^{-1}$.
\end{prop}

\begin{prop}\label{tubenbdprop}  Use the notation of Corollary \ref{tubenbdcor} and
Proposition \ref{Upsilonprop}.  There exist an open set
$U\subseteq\Lambda_+^2T^*N$ containing the zero section and
$W=(\jmath_N\circ\Upsilon)(V)$, a tubular neighbourhood $T$ of $N$
in $\R^7$ containing $S$, and a diffeomorphism $\delta:U\rightarrow
T$, affine on the fibres, that takes the zero section of\/
$\Lambda_+^2T^*N$ to $N$ and such that the following diagram
commutes:
\begin{equation}\label{deltaeq}
\begin{gathered}
\xymatrix{
W\ar[r]^{\tilde{\Upsilon}^{-1}}\ar[d]_{\delta} &
\jmath_C(V) \ar[d]^{\jmath_C^{-1}}
\\
S & V.\ar[l]_{n_P}
}\end{gathered}
\end{equation}
 Moreover, we may choose $U$ to be a $C^1_1$-open neighbourhood of the zero section. 
\end{prop}

\begin{proof} Define the diffeomorphism $\delta|_W:W\rightarrow S$ by
\eq{deltaeq}.  Interpolating smoothly over the compact set $K\subseteq N$, we extend $S$ to
$T$, $W$ to $U$ and $\delta|_W$ to $\delta$ as required.
Furthermore, by Corollary \ref{tubenbdcor}, $V$ can be chosen to be an open neighbourhood of the zero section in $C^1_1$.  
Hence, we can arrange the same for $U$.
\end{proof}

\subsection{The map {\boldmath $F$} and the associated map
{\boldmath $G$}}\label{FGsubsection}

We introduce the notation $C^{k}_{\lambda}(U)=\{\alpha\in
C^{k}_{\lambda}(\Lambda_+^2T^*N)\,:\,\Gamma_\alpha\subseteq U\},$ where $U$ is
given by Proposition \ref{tubenbdprop} and $\Gamma_{\alpha}$ is the graph of $\alpha$.  
The fact that $U$ is a $C^1_1$-open set  ensures that
$C^{k}_{\lambda}(U)$ is an open subset of
$C^{k}_{\lambda}(\Lambda_+^2T^*N)$ for $k\geq 1$, since $\lambda<1$.  We use
similar conventions to define subsets of the spaces discussed
in $\S$\ref{weightedsection}, though we must take care the spaces contain 
continuous forms so that their graphs are well-defined.  Moreover, these subsets will be open whenever the
space embeds continuously in $C^1_1$; for example, $L^p_{k+1,\,\lambda}(U)$ is an open neighbourhood of zero 
in $L^p_{k+1,\,\lambda}(\Lambda^2_+T^*N)$ if $p>4$ and $k\geq 1$ by Theorem \ref{wSobthm}.

We may now describe our \emph{deformation map}.

\begin{dfn}\label{Fdfn} Use the notation from the start of this section 
 and 
Proposition \ref{tubenbdprop}.  For $\alpha\in C^1_{\loc}( U)$ let
$\pi_{\alpha}:N\rightarrow\Gamma_{\alpha}\subseteq U$, where
$\Gamma_{\alpha}$ is the graph of $\alpha$, be
given by $\pi_{\alpha}(x)=\big(x,\alpha(x)\big)$. Let
$f_\alpha=\delta\circ\pi_{\alpha}$ and let
$N_{\alpha}=f_\alpha(N)$, so that
$N_{\alpha}$ is the deformation of $N$ corresponding to $\alpha$.
We 
 define
$F:C^{1}_{\loc}(U)\rightarrow
C^{0}_{\loc}(\Lambda^3T^*N)$ by
\begin{equation*}\label{ch7s1eq1} 
F(\alpha)=f_\alpha^*\big(\varphi|_{N_{\alpha}}\big). \end{equation*} By
Definition \ref{coassdfn}, ${\rm Ker}\, F$ is the set of $\alpha\in C^1_{\loc}(U)$ such that
the deformation $N_{\alpha}$ of $N$ is coassociative.  Note that $F$ is a nonlinear operator and that
$$dF|_0(\alpha)=d\alpha,$$ 
for all $\alpha\in C^1_{\loc}(\Lambda^2_+T^*N)$, by \cite[p.~731]{McLean} and our choice of $\delta$. 
\end{dfn}

However, we only wish to consider \emph{smooth} coassociative deformations
$N_{\alpha}$ of $N$ which are
\emph{asymptotically conical} to $C$ with rate $\lambda$.  If $N_{\alpha}$ is such a
deformation, then $\alpha\in C^{\infty}(U)$ and there exists a diffeomorphism
$\Psi_{\alpha}:(R,\infty)\times\Sigma \rightarrow
N_{\alpha}\setminus K_{\alpha}$, where $K_{\alpha}$ is a compact
subset of $N_{\alpha}$, as in Definition \ref{ACsubmflddfn}.  We may define $\Psi_{\alpha}$ such that
$\Psi_{\alpha}(r,\sigma)-\iota(r,\sigma)$ is orthogonal to
$T_{r\sigma}C$ for all $\sigma\in\Sigma$ and $r>R$.

Recall the notation of Corollary \ref{tubenbdcor} and Proposition \ref{Upsilonprop}.
Before Proposition \ref{Upsilonprop},  we
showed that $\Psi-\iota$ can be identified with the graph of
$\gamma_C\in C^{\infty}_\lambda(\Lambda^2_+T^*P)$ and therefore with
the graph of $\alpha_C=\tilde{\Upsilon}(\gamma_C)\in
C^\infty_\lambda(\Lambda^2_+T^*Q)$. Similarly, $\Psi_\alpha-\iota$
can be identified with the graph of
$\alpha+\alpha_C\in
C^\infty_\lambda(\Lambda^2_+T^*Q)$.  Hence, $\alpha\in
C^{\infty}_\lambda(U)$.

We conclude that $N_{\alpha}$ is AC to $C$ with rate $\lambda$ if
and only if $\alpha\in C_{\lambda}^{\infty}(U)$.  This leads immediately to our next result, which gives a local description
of the moduli space $\mathcal{M}(N,\lambda)$ using the deformation map $F$.

\begin{prop}\label{Fmodulispaceprop}
Using the notation of Definitions \ref{modulispacedfn} and
\ref{Fdfn}
, $\mathcal{M}(N,{\lambda})$ is locally homeomorphic
to the kernel of $F:C_{\lambda}^{\infty}(U)\rightarrow
C^{\infty}(\Lambda^3T^*N)$.
\end{prop}

We now prove 
 that we can rewrite $F(\alpha)$ as a sum of $d\alpha$ and a term which is no worse than quadratic in $\alpha$
and $\nabla\alpha$.  This result will be useful throughout the article, but in particular we shall need it to 
derive regularity results in $\S$\ref{regsubsection}.  

\begin{prop}\label{PFprop} Use the notation of Definitions \ref{wSobdfn}-\ref{wHolderdfn} and \ref{Fdfn}. 
We can write
\begin{equation}\label{PFeq}
F(\alpha)(x)=d\alpha(x)+P_F\big(x,\alpha(x),\nabla\alpha(x)\big)
\end{equation}
for $x\in N$, where 
$$P_F:\{(x,y,z)\,:\,(x,y)\in U,\,z\in
T_x^*N\otimes\Lambda_+^2T_x^*N\}\rightarrow\Lambda^3T^*N$$ 
is a
smooth map such that $P_F(x,y,z)\in\Lambda^3T^*_xN$.  Denote $P_F\big(x,\alpha(x),\nabla\alpha(x)\big)$ by
$P_F(\alpha)(x)$ for all $x\in N$.  

Let $\mu<1$.  For each $n\in\N$, if $\alpha\in C^{n+1}_\mu(U)$ and
$\|\alpha\|_{C^1_1}$ is sufficiently small, $P_F(\alpha)\in
C^n_{2\mu-2}(\Lambda^3T^*N)$ and there exists a constant $c_n>0$
such that
\begin{equation}\label{PFestckeq}
\|P_F(\alpha)\|_{C_{2\mu-2}^n}\leq
c_n\|\alpha\|_{C^{n+1}_{\mu}}^2\,.
\end{equation}
Let $p>4$, $k\geq 1$, $l\in\N$ and $a\in(0,1)$.  
If $\alpha\in L^p_{k+1,\,\mu}(U)$ or $\alpha\in C^{l+1,\,a}_{\mu}(U)$,
with $\|\alpha\|_{C^1_1}$ sufficiently small, $P_F(\alpha)\in L^p_{k,\,2\mu-2}(\Lambda^3T^*N)$ or 
$C^{l,\,a}_{2\mu-2}(\Lambda^3T^*N)$ and there exist constants $l_{p,\,k}>0$ and $c_{l,\,a}>0$ such that 
$$\|P_F(\alpha)\|_{L^p_{k,\,2\mu-2}}\leq l_{p,\,k}\|\alpha\|_{L^p_{k+1,\,\mu}}^2
\quad\text{or}\quad \|P_F(\alpha)\|_{C^{l,\,a}_{2\mu-2}}\leq c_{l,\,a}\|\alpha\|_{C^{l+1,\,a}_{\mu}}^2.$$
\end{prop}

\begin{remarks}
As $\mu<1$, $2\mu-2<\mu-1$, so $C^n_{2\mu-2}\hookrightarrow C^n_{\mu-1}$.  Similar  
continuous embeddings occur for the weighted Sobolev and H\"older spaces.  Furthermore, the conditions 
$p>4$, $k\geq 1$ and $l\in\N$ ensure, by Definition \ref{wHolderdfn} and Theorem \ref{wSobthm}(b), that 
$C^{l+1,\,a}_{\mu}\hookrightarrow C^1_1$ and $L^p_{k+1,\,\mu}\hookrightarrow C^1_1$.
\end{remarks}

\begin{proof} First, by the definition of $F$, $F(\alpha)(x)$
relates to the tangent space to the graph $\Gamma_{\alpha}$ of $\alpha$ at
$\pi_{\alpha}(x)=\big(x,\alpha(x)\big)$.  Note that $T_{\pi_{\alpha}(x)}\Gamma_{\alpha}$
depends on both $\alpha(x)$ and $\nabla\alpha(x)$ and hence so must
$F(\alpha)(x)$.  We may then define $P_F$ by \eq{PFeq} such that
it is a smooth function of its arguments as claimed.

We only prove the esimate \eq{PFestckeq} on $P_F$ as the results for weighted Sobolev and H\"older spaces can be deduced from
the work presented here.
Recall the notation at the start of this section and of Corollary \ref{tubenbdcor} and Proposition \ref{Upsilonprop}.  Let
$\nabla_{C}$ denote the Levi--Civita connection of the conical metric on $C$.  
We argued after Corollary \ref{tubenbdcor} that we may identify the
displacement $\Psi-\iota$ of $N$ from $C$, outside the
compact subset $K$, with $\gamma_C\in
C_{\mu}^{\infty}(\Lambda_+^2T^*P)$.  For each $\alpha\in C^1_{\loc}(U)$, there exists
a unique $\gamma\in C^{1}_{\loc}(\Lambda^2_+T^*P)$ such that $\alpha=\tilde{\Upsilon}(\gamma)$ on
$N\setminus K$.  Thus, define a function
$F_C(\gamma+\gamma_C)$, for $\gamma\in
C_{\loc}^{1}(\Lambda_+^2T^*P)$, on $(R,\infty)\times\Sigma$
by
\begin{equation}\label{FCeq}
F_C(\gamma+\gamma_C)(r,\sigma)=F(\alpha)
\big(\Psi(r,\sigma)\big),
\end{equation}
where $\alpha|_{N\setminus K}=\tilde{\Upsilon}(\gamma)$.
Now define a smooth function $P_C$ by an equation analogous to
\eq{PFeq}:
\begin{align}
F_C(\gamma&+\gamma_C)(r,\sigma)\nonumber
\\
&=
d(\gamma+\gamma_C)(r,\sigma)+P_C\big(
(r,\sigma),(\gamma+\gamma_C)(r,\sigma),\nabla_C(\gamma+\gamma_C)(r,\sigma)\big).\label{PCeq}
\end{align}
We notice that $F_C$ and $P_C$ are only dependent on the cone $C$
and, rather trivially, on $R$.  Therefore, because of this fact and
our choice of $\delta$ in Proposition \ref{tubenbdprop}, these
functions have scale equivariance properties.  We may therefore
derive equations and inequalities on $\{R\}\times\Sigma$ and deduce
the result on all of $(R,\infty)\times\Sigma$ by introducing an
appropriate scaling factor of $r$.

Now, since the graph of $\alpha=0$ corresponds to our coassociative 4-fold $N$,
$F\big(\tilde\Upsilon(0)\big)=F(0)=0$.  So, taking $\gamma=0$ in \eq{FCeq} gives:
\begin{equation}\label{FCPCeq}
F_C(\gamma_C)=d\gamma_C+P_C(\gamma_C)=0, \end{equation} adopting
similar notation for $P_C(\gamma_C)$ as for
$P_F(\alpha)$.  From \eq{PFeq}-\eq{FCPCeq} we calculate:
\begin{align}
P_F(\alpha)\big(\Psi(r,\sigma)\big)&
=d\gamma_C(r,\sigma)+P_C(\gamma+\gamma_C)(r,\sigma)\nonumber\\ 
&=d\gamma_C(r,\sigma)+P_C(\gamma+\gamma_C)(r,\sigma)-\big(d\gamma_C(r,\sigma)+P_C(\gamma_C)(r,\sigma)\big) \nonumber\\
\label{PFPCeq}
&=P_C(\gamma+\gamma_C)(r,\sigma)-P_C(\gamma_C)(r,\sigma).
\end{align}
Noting that $P_C$ is a function of three variables $x$, $y$ and
$z$, we see that
\begin{align}
P_C(\gamma+\gamma_C)-P_C(\gamma_C)&=\int_0^1\frac{d}{dt}\,P_C(t\gamma+\gamma_C)
\,dt
\nonumber\\[4pt]
\label{PCest1eq}
 & =\int_0^1\!\gamma\cdot\frac{\partial P_C}{\partial
y}(t\gamma+\gamma_C)+\nabla_C\gamma\cdot\frac{\partial P_C}{\partial
z}(t\gamma+\gamma_C)\,dt.
\end{align}
By Taylor's Theorem,
\begin{equation}\label{PCest2eq}
P_C(\gamma+\gamma_C)=P_C(\gamma_C)+\gamma\cdot\frac{\partial
P_C}{\partial y}(\gamma_C)+ \nabla_C\gamma\cdot\frac{\partial
P_C}{\partial z}(\gamma_C)
+O(r^{-2}|\gamma|^2+|\nabla_C\gamma|^2)\end{equation} when $r^{-1}|\gamma|$
and $|\nabla_C\gamma|$ are small. 

Since
$dF|_0(\alpha)=d\alpha$, as noted in
Definition \ref{Fdfn},
$dF_C|_{\gamma_C}(\gamma+\gamma_C)=d\gamma$ and hence
$dP_C|_{\gamma_C}=0$.  Thus, the first derivatives of $P_C$ with
respect to $y$ and $z$ must vanish at $\gamma_C$ by \eq{PCest2eq}.
Therefore, given small
$\epsilon>0$, there exists a constant $A_0>0$ such that
\begin{align}
\label{PCest3aeq}
\left|\frac{\partial P_C}{\partial y}(t\gamma+\gamma_C)\right| &\leq
A_0(r^{-2}|\gamma| +r^{-1}|\nabla_C\gamma|)\quad\text{and}\\[4pt]
\label{PCest3beq}
\left|\frac{\partial P_C}{\partial z}(t\gamma+\gamma_C)\right| &\leq
A_0(r^{-1}|\gamma| +|\nabla_C\gamma|)
\end{align} for $t\in [0,1]$,
whenever
\begin{equation}\label{epsiloneq}
r^{-1}|\gamma|,\,r^{-1}|\gamma_C|,\,|\nabla_C\gamma|\;\text{and}\;|\nabla_C\gamma_C|\leq\epsilon.
\end{equation}
The factors of $r$ are determined by considering the scaling properties of $\gamma$ and $P_C$ and their derivatives under
changes in $r$.

By \eq{gammaCeq}, $r^{-1}|\gamma_C|$ and $|\nabla_C\gamma_C|$ tend to
zero as $r\rightarrow\infty$.  We can thus ensure that \eq{epsiloneq}
is satisfied by the $\gamma_C$ terms by making $R$ larger.
Hence, \eq{epsiloneq} holds if
$\|\gamma\|_{C^1_1}\leq\epsilon$.
 Therefore, putting estimates
\eq{PCest3aeq} and \eq{PCest3beq} in \eq{PCest1eq} and using \eq{PFPCeq},
\begin{align}
\big|P_F(\alpha)\big(\Psi(r,\sigma)\big)\big|&=|P_C(\gamma+\gamma_C)(r,\sigma)-P_C(\gamma_C)(r,\sigma)|\nonumber\\
&\leq
A_0\big(r^{-1}|\gamma(r,\sigma)|+|\nabla_C\gamma(r,\sigma)|\big)^2\label{PFest1eq}
\end{align}
for all $r>R$, $\sigma\in\Sigma$,
whenever $\|\gamma\|_{C^1_1}\leq\epsilon$.   

Therefore, if $\alpha\in C^{1}_{\mu}(U)$, the corresponding 
$\gamma$ lies in $C^{1}_{\mu}(\Lambda^2_+T^*P)$.  Hence, by \eq{PFest1eq}, $P_F(\alpha)\in C^0_{\mu}$ and 
there exists a constant $c_0$ such that
$$\sup_N|\rho^{2-2\mu}P_F(\alpha)|\leq c_0 \left(\sum_{i=0}^1\sup_N|\rho^{i-\mu}\nabla^i\alpha|\right)^{\!2}$$
whenever $\|\alpha\|_{C^1_1}$ is sufficiently small, where $\rho$ is a radius function on $N$.  Thus, \eq{PFestckeq} holds for $n=0$.

Similar calculations give analogous results to \eq{PFest1eq} for
derivatives of $P_F$, from which we can deduce \eq{PFestckeq} for $n>0$.   We shall explain the method by considering
the first derivative. If $\gamma\in C^2_{\loc}(\Lambda^2_+T^*P)$, we calculate from \eq{PCest1eq}:
\begin{align*}
\nabla_C&\big(P_C(\gamma+\gamma_C)-P_C(\gamma_C)\big)\\[4pt]
=&\int_0^1
\nabla_C\!\left(\gamma\cdot\frac{\partial P_C}{\partial
y}(t\gamma+\gamma_C)+\nabla_C\gamma\cdot\frac{\partial
P_C}{\partial z}(t\gamma+\gamma_C)\right)dt \\[4pt]
=&\int_0^1\nabla_C\gamma\cdot\frac{\partial P_C}{\partial
y}+\gamma\cdot\!\left(\nabla_C(t\gamma+\gamma_C)\cdot\frac{\partial^2
P_C}{\partial y^2}+\nabla^2_C(t\gamma+\gamma_C)\cdot\frac{\partial^2
P_C}{\partial y\partial z}\right) \\[4pt]
&+\nabla^2_C\gamma\cdot\frac{\partial P_C}{\partial
z}+\nabla_C\gamma\cdot\!\left(\nabla_C(t\gamma+\gamma_C)\cdot\frac{\partial^2
P_C}{\partial z\partial
y}+\nabla^2_C(t\gamma+\gamma_C)\cdot\frac{\partial^2 P_C}{\partial
z^2}\right)dt.
\end{align*}
Whenever $\|\gamma\|_{C^1_1}\leq\epsilon$ there exists a constant
$A_1>0$ such that \eq{PCest3aeq} and \eq{PCest3beq} hold with $A_0$ replaced by $A_1$
and, for $t\in[0,1]$,
$$\left|\frac{\partial^2
P_C}{\partial y^2}(t\gamma+\gamma_C)\right|,\,\left|\frac{\partial^2
P_C}{\partial y\partial z}(t\gamma+\gamma_C)\right|\;\text{and}\;
\left|\frac{\partial^2 P_C}{\partial
z^2}(t\gamma+\gamma_C)\right|\leq A_1,$$ since the second
derivatives of $P_C$ are continuous functions defined on the closed
bounded set given by $\|\gamma\|_{C^1_1}\leq\epsilon$.  
We see that
\begin{align*}
\big|\nabla
\big(P_F(\alpha)\big)\big(\Psi(r,\sigma)\big)\big|&
=\big|\nabla_C\big(P_C(\gamma+\gamma_C)-P_C(\gamma_C)\big)(r,\sigma)\big|\\[4pt]
&\leq
A_1r\left(\sum_{i=0}^2r^{i-2}\big|\nabla_C^i\gamma(r,\sigma)\big| \right)^{\!2}
\end{align*}
whenever $\|\gamma\|_{C^1_1}\leq\epsilon$.  From this we can deduce the result \eq{PFestckeq} for $n=1$.

In general we have the estimate
$$\big|\nabla^j\big(P_F(\alpha)\big)\big(\Psi(r,\sigma)\big)\big|\leq A_j
r^j\left(\sum_{i=0}^{j+1}r^{i-(j+1)}\big|\nabla^i_C\gamma(r,\sigma)\big|\right)^{\!2}$$ for
some $A_j>0$ whenever $\|\gamma\|_{C^1_1}\leq\epsilon$.  The result \eq{PFestckeq} for all $n\in\N$
follows.
\end{proof}

We now turn to the associated map $G$.

\begin{dfn}\label{Gdfn} Use the notation of 
Proposition \ref{tubenbdprop} and Definition \ref{Fdfn}.  
Define $G:C_{\loc}^{1}(U)\times
C_{\loc}^{1}(\Lambda^4T^*N)\rightarrow
C^{0}_{\loc}(\Lambda^3T^*N)$ by:
\begin{equation*}
G(\alpha,\beta)=F(\alpha)+d^*\beta. 
\end{equation*} 
By the observations made in Definition \ref{Fdfn},
\begin{align*}
dG|_{(0,0)}
(\alpha,\beta)= d\alpha+d^*\beta
\end{align*}
for all $(\alpha,\beta)\in C^1_{\loc}(\Lambda^2_+T^*N\oplus\Lambda^4T^*N)$.  Thus $G$ is a 
nonlinear elliptic operator at $(0,0)$; that is, the linearisation of $G$ at $(0,0)$ is elliptic.
\end{dfn}

To complete this subsection, we relate the kernels of $F$ and $G$.

\begin{prop}\label{isokernelsprop}  Use the notation of Definitions \ref{wSobdfn}-\ref{wHolderdfn}, \ref{Fdfn} and \ref{Gdfn}.  
Let $p>4$, $k\geq 1$, $l\in\N$, $a\in(0,1)$, $\mu<0$ and $\nu\leq 0$.  
The kernels of $F$ and $G$ in $C^{l+1}_{\mu}$, $C^{l+1,\,a}_{\mu}$ and $L^p_{k+1,\,\nu}$ respectively
 are isomorphic.
\end{prop}

\begin{proof} 
Note that $\varphi$ is exact on any deformation $N_{\alpha}$ of $N$
because it is closed and vanishes on $N$. 
Thus, if $G(\alpha,\beta)=0$,
$$d\big(G(\alpha,\beta)\big)=d\big(F(\alpha)\big)+dd^*\beta=\Delta\beta=0,$$
since $F(\alpha)$ is exact.  Recall we have a radius function $\rho$ on $N$. 
If $\beta$ decays with order $O(\rho^\mu)$, or $o(\rho^{\nu})$, as
$\rho\rightarrow\infty$, then $*\beta$ is a
harmonic function on $N$ which tends to zero as
$\rho\rightarrow\infty$.  The Maximum Principle allows us to deduce that $*\beta=0$ and conclude
that $\beta=0$, from which the proposition follows.
\end{proof}

\subsection{Uniformly elliptic AC operators and regularity}\label{regsubsection}

We want to consider the regularity of solutions $F(\alpha)=0$ near $0$, but this equation is not elliptic.  Therefore, we study
 $G(\alpha,\beta)=0$ near $(0,0)$, which is a \emph{nonlinear}
elliptic equation.  To do so, we must first discuss the regularity theory of certain \emph{linear} elliptic operators,
which we shall define, acting between
weighted Banach spaces.  We shall use the theory in \cite{Lockhart} and \cite{LockhartMcOwen}, which centres around \emph{asymptotically 
cylindrical} manifolds, so we begin with their definition.

\begin{dfn}\label{ACyldfn}
Recall Definition \ref{wSobdfn} and the notation introduced at the start of $\S$\ref{defmapsection}.  In particular, 
we have a compact subset $K$ of $N$ such that $N_{\infty}=N\setminus K\cong (R,\infty)\times\Sigma$ 
and a radius function $\rho$ on $N$.  
Let $g$ be the metric on $N$,
considered as a Riemannian manifold.  Define a metric $\tilde{g}$ on $N$ by $\tilde{g}=\rho^{-2}g$.  
Further, define a coordinate $t$ on $(T,\infty)$, where $T=\log R$, by
$t=\log r$, let $g_{\text{cyl}}=dt^2+g_{\Sigma}$ be the \emph{cylindrical} metric on 
$(T,\infty)\times\Sigma$ and let $\tilde{\nabla}$ be the Levi--Civita connection of 
$\tilde{g}$.  We say that $(N,\tilde{g})$
is \emph{asymptotically cylindrical} (ACyl).

For $p\geq 1$, $k\in\N$ and $\mu\in\R$ we define the Banach space $\tilde{L}^p_{k,\,\mu}(\Lambda^mT^*N)$ to be the subspace of 
$L^p_{k,\,\loc} (\Lambda^mT^*N)$ such that the following norm is finite:
$$\|\xi\|_{\tilde{L}^p_{k\,\,\mu}}=\left(\sum_{j=0}^k\int_{N}|\rho^{-\mu}\tilde{\nabla}^j\xi|^p
dV_{\tilde{g}}\right)^{\frac{1}{p}}.$$
\end{dfn}

\begin{dfn}\label{ACylopdfn} Use the notation from 
Definition \ref{ACyldfn}.  
Suppose that $\mathcal{P}:C^l_{\loc}(\Lambda^mT^*N)\rightarrow C^0_{\loc}(\Lambda^{m^{\prime}}T^*N)$ 
 and $\mathcal{P}_{\infty}:C^l_{\loc}(\Lambda^mT^*N_{\infty})\rightarrow C^0_{\loc}(\Lambda^{m^{\prime}}T^*N_{\infty})$ are
 linear differential operators of order $l$.  Suppose further that $\mathcal{P}_{\infty}$
is invariant under the $\R^+$-action on $N_{\infty}\cong (T,\infty)\times\Sigma$.  
We say that $\mathcal{P}_{\infty}$ is \emph{cylindrical}.  

For $\xi\in C^l_{\loc}(\Lambda^mT^*N_{\infty})$, 
$$\mathcal{P}\xi=\sum_{i=0}^l\mathcal{P}_i\cdot\tilde{\nabla}^i\xi\quad\text{and}\quad\mathcal{P}_{\infty}\xi=\sum_{i=0}^l
\mathcal{P}_{i,\,\infty}\cdot\tilde{\nabla}^i\xi,$$
where $\mathcal{P}_i$ and $\mathcal{P}_{i,\,\infty}$ are tensors on $N_{\infty}$ of type $(m+i,m^{\prime})$ and 
``$\,\cdot\,$'' means tensor product followed by contraction.  
If, for $i=0,\ldots, l$,
$$|\tilde{\nabla}^j(\mathcal{P}_i-\mathcal{P}_{i,\,\infty})|\rightarrow 0\qquad \text{for $j\in\N$ as $t\rightarrow \infty$,}$$
where $|\,.\,|$ is calculated using $g_{\text{cyl}}$,
we say that $\mathcal{P}$ is \emph{asymptotically cylindrical} (to $\mathcal{P}_{\infty}$).
By \cite[Theorem 3.7]{Lockhart}, $\mathcal{P}$ is a continuous map from $\tilde{L}^p_{k+l,\,\mu}(\Lambda^mT^*N)$ to $\tilde{L}^p_{k,\,\mu}(\Lambda^{m^{\prime}}T^*N)$
for all $p>1$, $k\in\N$ and $\mu\in\R$.
\end{dfn}

\begin{dfn}\label{ACopdfn}
Use the notation from Definitions \ref{ACyldfn} and \ref{ACylopdfn}, so we have an operator $\mathcal{P}$ of 
order $l$ acting between
$m$- and $m^{\prime}$-forms on $N$.  Let $\nu\in\R$.  
We say that $\mathcal{P}$ is \emph{asymptotically conical} (AC) with \emph{rate} $\nu$ if the differential operator 
$$\mathcal{P}^{\nu}=\rho^{-m^{\prime}+\nu}\mathcal{P}\rho^{m}$$
is asymptotically cylindrical to a cylindrical map $\mathcal{P}_{\infty}$, say.  By 
\cite[Proposition and Definition 4.4]{Lockhart}, $\mathcal{P}:L^p_{k+l,\,\mu}(\Lambda^mT^*N)\rightarrow L^p_{k,\,\mu-\nu}
(\Lambda^{m^{\prime}}T^*N)$ is continuous for all $p>1$, $k\in\N$ and $\mu\in\R$. 
\end{dfn}

\begin{note}
This definition can clearly be extended to linear differential operators acting between more general bundles of forms.
\end{note}

\noindent Notice that an AC operator with rate $\nu$ reduces the
growth rate of a form on the ends of $N$ by $\nu$. Examples of AC operators abound: $d$ and $d^*$ are first-order operators with
rate $1$ and the Laplacian is a second-order operator with rate $2$.

\begin{dfn}\label{uniformdfn}  Use the notation of Definition
\ref{ACopdfn}.  We say that an AC operator $\mathcal{P}$ is
\emph{uniformly elliptic} if it is elliptic and $\mathcal{P}_{\infty}$ is
elliptic.
\end{dfn}

\begin{remark} The definition of uniformly elliptic above
implies uniform ellipticity in the sense of global bounds on the
coefficients of the symbol.
\end{remark}  

\noindent The operators $d+d^*$ and the Laplacian
are 
 uniformly elliptic AC operators.

\smallskip

We now turn to regularity results for \emph{smooth} uniformly elliptic AC operators.

\begin{thm}\label{ACopregthm} Let $V$ and $W$ be bundles of forms on
$N$ and let $\mathcal{P}$ be a smooth uniformly elliptic AC operator 
from $V$ to $W$ of order $l$ and rate $\nu$, in the sense of Definitions \ref{ACopdfn} and \ref{uniformdfn}.  Let 
$p>1$, $k\in\N$, $a\in(0,1)$ and $\mu\in\R$.

\begin{itemize}\item[\emph{(a)}] 
Suppose that $\mathcal{P}\xi=\eta$ holds for $\xi\in
L^1_{l,\,\loc}(V)$ and $\eta\in
L^1_{0,\,\loc}(W)$. If\/ $\xi\in L^p_{0,\,\mu}(V)$ and
$\eta\in L^p_{k,\,\mu-\nu}(W)$, then $\xi\in L^p_{k+l,\,\mu}(V)$ and
$$\|\xi\|_{L^p_{k+l,\,\mu}}\leq c\left(\|\eta\|_{L^p_{k,\,\mu-\nu}}+
\|\xi\|_{L^p_{0,\,\mu}}\right)$$ for some constant $c>0$ independent
of $\xi$ and $\eta$.

\item[\emph{(b)}]
Suppose that $\mathcal{P}\xi=\eta$ holds for $\xi\in
C^l_{\loc}(V)$ and $\eta\in
C^0_{\loc}(W)$. If\/ $\xi\in C^0_{\mu}(V)$ and $\eta\in
C^{k,\,a}_{\mu-\nu}(W)$, then $\xi\in C^{k+l,\,a}_{\mu}(V)$ and
$$\|\xi\|_{C^{k+l,\,a}_{\mu}}\leq c^\prime\left(\|\eta\|_{C^{k,\,a}_{\mu-\nu}}
+\|\xi\|_{C^0_{\mu}}\right)$$ for some constant $c^\prime>0$
independent of $\xi$ and $\eta$.  Moreover, these estimates hold if
the coefficients of $\mathcal{P}$ only lie in $C^{k,\,a}_{\loc}$.
\end{itemize}
\end{thm}

\noindent These results may be deduced from those given in \cite[$\S$6.1.1]{Marshall} or \cite{Mazya}.

By taking $\eta=0$ in Theorem \ref{ACopregthm}(b), we have the following useful corollary.

\begin{cor}\label{ACopkerregcor} In the notation of Theorem
\ref{ACopregthm}, if\/ $\xi\in C^l_{\mu}(V)$ satisfies $\mathcal{P}\xi=0$ then
$\xi\in C^{\infty}_{\mu}(V)$.
\end{cor}

We conclude by achieving the aim of this subsection.

\begin{prop}\label{Gregprop} $\!$Use the notation of Definition \ref{Gdfn}. 
Let $(\alpha,\beta)\in
L_{k+1,\,\mu}^p(U)\times L_{k+1,\,\mu}^p(\Lambda^4T^*N)$ for 
$p>4$, $k\geq 2$ and $\mu<1$.  If $G(\alpha,\beta)=0$ and
$\|\alpha\|_{C^1_1}$ is sufficiently small, $(\alpha,\beta)\in
C_{\mu}^{\infty}(U)\times C_{\mu}^{\infty}(\Lambda^4T^*N)$.
\end{prop}

\begin{proof}
 Notice first that 
 $\alpha$ and $\beta$ lie in $C^2_{\mu}$ by Theorem
\ref{wSobthm}, since $k+1-\frac{4}{p}
>2$.  

As noted in the proof of Proposition \ref{isokernelsprop}, 
$G(\alpha,\beta)=0$ implies that $\Delta\beta=0$.  Hence, by 
Corollary \ref{ACopkerregcor}, $\beta\in
C^{\infty}_\mu(\Lambda^4T^*N)$.



For the following argument we find it useful to work with weighted
H\"older spaces, defined in Definition \ref{wHolderdfn}.  By Theorem \ref{wSobthm}, $\alpha\in
C^{k,\,a}_{\mu}(U)$ with $a=1-4/p\in (0,1)$ since $p>4$.  Let $\pi_{\Lambda^2_+}$ be the projection from 
2-forms to self-dual 2-forms on $N$.  We 
know that $d^*\big(G(\alpha,\beta)\big)=d^*\big(F(\alpha)\big)=0$ 
and also that 
$$\tilde{F}(\alpha)=\pi_{\Lambda^2_+}\Big(d^*\big(F(\alpha)\big)\Big)=0$$
is a nonlinear \emph{elliptic} equation at $0$
, meaning that $d\tilde{F}|_0$
 is an elliptic operator. 

We can write $\tilde{F}$ as
$$\tilde{F}(\alpha)(x)=R_F\big(x,\alpha(x),\nabla\alpha(x)\big)\nabla^2\alpha(x)
+E_F\big(x,\alpha(x),\nabla\alpha(x)\big),$$ where
$R_F
$ and
$E_F
$ are smooth functions of their
arguments, since $\tilde{F}(\alpha)$ is \emph{linear} in $\nabla^2\alpha$ with coefficents depending on $\alpha$ and $\nabla\alpha$.  
Define
$$S_{\alpha}(\gamma)(x)= R_F\big(x,\alpha(x),\nabla\alpha(x)\big)
\nabla^2\gamma(x)$$ for $\gamma\in
C^2_{\loc}(\Lambda^2_+T^*N)$.  Note that $S_{\alpha}$ is \emph{not} the linearisation of $\tilde{F}$. 
Then $S_{\alpha}$ is a 
\emph{linear uniformly elliptic} second order AC operator with rate 2, if
$\|\alpha\|_{C^1_1}$ is sufficiently small, whose coefficients
depend on $x$, $\alpha(x)$ and $\nabla\alpha(x)$. These coefficients
therefore lie in $C^{k-1,\,a}_{\loc}$.  

Recall the notation and
results of Proposition \ref{PFprop} and that $d^*$ is an AC operator with rate 1 in the sense of
Definition \ref{ACopdfn}.  Thus, 
$d^*d\alpha+d^*\big(P_F(\alpha)\big)=0$ and $d^*\big(P_F(\alpha)\big)\in
C^{k-2,\,a}_{2\mu-3}(\Lambda^2T^*N)$.  
Therefore,
$$S_{\alpha}(\alpha)(x)=-E_F\big(x,\alpha(x),\nabla\alpha(x)\big)\in
C^{k-2,\,a}_{2\mu-3}(\Lambda^2_+T^*N)\subseteq
C^{k-2,\,a}_{\mu-2}(\Lambda^2_+T^*N),$$ since $\mu<1$.  
However, $E_F\big(x,\alpha(x),\nabla\alpha(x)\big)$ only depends on $\alpha$
and $\nabla\alpha$, and is at worst quadratic in these quantities by
Proposition \ref{PFprop}, so it must in fact lie in
$C^{k-1,\,a}_{\mu-2}(\Lambda^2_+T^*N)$ since we are given control
on the decay of the first $k$ derivatives of $\alpha$ at infinity.

As $\alpha\in C^2_{\mu}(\Lambda_+^2T^*N)$ and
$S_{\alpha}(\alpha)\in C^{k-1,\,a}_{\mu-2}(\Lambda^2_+T^*N)$,
Theorem \ref{ACopregthm}(b) implies that $\alpha\in
C^{k+1,\,a}_{\mu}(\Lambda_+^2T^*N)$.  Therefore we have proved $\alpha\in
C^{k+1,\,a}_{\mu}(\Lambda^2_+T^*N)$ only knowing a priori that
$\alpha\in C^{k,\,a}_{\mu}(\Lambda^2_+T^*N)$.  We proceed by
induction to show that $\alpha\in
C^{l,\,a}_{\mu}(\Lambda^2_+T^*N)$ for all $l\geq 2$.  The result follows from the elementary
observation made at the end of Definition \ref{wHolderdfn}.
\end{proof}

Taking $\beta=0$ in Proposition \ref{Gregprop} gives our main regularity result.
\begin{cor}\label{Fregcor} Use the notation of Definition \ref{Fdfn}. 
Let $\alpha\in
L_{k+1,\,\mu}^p(U)$ for
$p>4$, $k\geq 2$ and $\mu<1$.  If $F(\alpha)=0$ and
$\|\alpha\|_{C^1_1}$ is sufficiently small, $\alpha\in
C_{\mu}^{\infty}(U)$.
\end{cor}

\section{The map {\boldmath $d+d^*$} and the exceptional set {\boldmath $\mathcal{D}$}}
\label{dd*section}


We begin by defining the map of interest. 

\begin{dfn}\label{dd*dfn}  
Let $p>4$, $k\geq 2$ and $\mu<1$.
Define the linear elliptic operator
\begin{equation}\label{dd*eq}
(d+d^*)_{\mu}:L^p_{k+1,\,\mu}(\Lambda^2_+T^*N\oplus\Lambda^4T^*N)\rightarrow L^p_{k,\,\mu-1}(\Lambda^3T^*N)
\end{equation}
by $(d+d^*)_{\mu}(\alpha,\beta)=d\alpha+d^*\beta$.  Let $\mathcal{K}(\mu)=\Ker(d+d^*)_{\mu}$.
\end{dfn}

\begin{remark}
The operator $(d+d^*)_{\mu}$ acts between the weighted Banach spaces claimed because it is AC with rate 1.
\end{remark}

We now make an important observation.

\begin{lem}\label{Gmappinglem}  Use the notation from Proposition \ref{tubenbdprop} 
and Definitions \ref{Gdfn} and \ref{dd*dfn}.
 By making the open set $U$ smaller in $C^1_1$ if necessary, 
\begin{equation}\label{Gmappingeq}
G:L_{k+1,\,\mu}^p(U)\times
L_{k+1,\,\mu}^p(\Lambda^4T^*N)\rightarrow
L_{k,\,\mu}^p(\Lambda^3T^*N).
\end{equation}
Moreover, the
linearisation of \eq{Gmappingeq} at $(0,0)$ acts as $(d+d^*)_{\mu}$, as in \eq{dd*eq}.
\end{lem}

\begin{proof}
By Proposition \ref{PFprop}, and the fact that $d$ and $d^*$ are AC operators with rate 1 on $N$, 
 we see that 
$G$ maps $(\alpha,\beta)\in L^p_{k+1,\,\mu}(U)\times L^p_{k+1,\,\mu}(\Lambda^4T^*N)$ into
$L^p_{k,\,\mu-1}(\Lambda^3T^*N)$ if $\|\alpha\|_{C^1_1}$ is sufficently small.  This bound on the norm of $\alpha$ can
be ensured by making the $C^1_1$-open set $U$ smaller.  The description of the linearisation follows from 
the observations in Definition \ref{Gdfn}.
\end{proof}

It is clear from Propositions \ref{Fmodulispaceprop} and \ref{isokernelsprop} that the kernel of the nonlinear map \eq{Gmappingeq}, 
for $\mu=\lambda$, is
 intimately linked with the moduli space $\mathcal{M}(N,\lambda)$.  Therefore, locally, 
one would expect the kernel $\mathcal{K}(\lambda)$ of the linear map  $(d+d^*)_{\lambda}$, given by
\eq{dd*eq}, to be related to $\mathcal{M}(N,\lambda)$ as well by Lemma \ref{Gmappinglem}.  
In fact, we shall see that $\mathcal{K}(\lambda)$ 
is directly connected with the \emph{infinitesimal deformations} of $N$.

\subsection{Fredholm and index theory}

We want to understand the Fredholm theory of \eq{dd*eq} and so state a result adapted from \cite[Theorem 1.1 \&
Theorem 6.1]{LockhartMcOwen}.
\begin{thm}\label{Pexceptthm}
Let $V$ and $W$ be bundles of forms over $N$,
let $p>1$, let $\mu,\nu\in\R$ and let $k,l\in\N$. Let
$\mathcal{P}:L^p_{k+l,\,\mu}(V)\rightarrow L^p_{k,\,\mu-\nu}(W)$ be a
uniformly elliptic AC operator of order $l$ and
rate $\nu$ in the sense of Definitions \ref{ACopdfn}-\ref{uniformdfn}. There exists a countable discrete set
$\mathcal{D}(\mathcal{P})\subseteq\R$, depending only on $\mathcal{P}_{\infty}$ as in
Definition \ref{ACopdfn}, such that $\mathcal{P}$ is Fredholm if and only if
$\mu\notin\mathcal{D}(\mathcal{P})$.
\end{thm}

From this we know that, for
each $m\in\N$ with $m\leq4$, there exists a countable discrete
subset $\mathcal{D}(\Delta^m)$ of $\R$ such that the Laplacian on
$m$-forms
$$\Delta^m:L_{k+2,\,\mu+1}^{p}(\Lambda^mT^*N)\rightarrow
L_{k,\,\mu-1}^{p}(\Lambda^mT^*N)$$ is Fredholm if and only if
$\mu+1\notin\mathcal{D}(\Delta^m)$.  
 Thus, $(d+d^*)_{\mu}$ is Fredholm if $\mu\notin\big(\mathcal{D}(\Delta^2)\cup\mathcal{D}(\Delta^4)\big)$.
However, we can give an explicit description of the set
$\mathcal{D}$ for which \eq{dd*eq} is not Fredholm,
following \cite[$\S$6.1.2]{Marshall}.

\medskip

Recall the notation from Definitions \ref{ACyldfn}-\ref{ACopdfn} and 
consider the maps $d$ and $d^*$ acting on $m$-forms on $N$.  
These are asymptotically conical with rate 1, so 
the related asymptotically cylindrical operators $d^1$ and $(d^*)^1$ are given by
$$d^1=\rho^{-m}d\rho^m\quad\text{and}\quad (d^*)^1=\rho^{-m+2}d^*\rho^m.$$
Notice that $\rho$ is asymptotic to $r=e^t$, where the cylindrical coordinates $(t,\sigma)$ on the
ends $N_{\infty}\cong (T,\infty)\times\Sigma$ of $N$ were introduced in 
Definition \ref{ACyldfn}.  Thus, we can take the cylindrical operator $(d+d^*)_{\infty}$ associated to $d+d^*$ 
to be
\begin{equation}\label{dd*inftyeq}
(d+d^*)_{\infty}=e^{-mt}(d+e^{2t}d^*)e^{mt},
\end{equation}
acting on $m$-forms on the ends of $N$.  

Since 
$N_{\infty}\cong (T,\infty)\times\Sigma$, an $m$-form $\alpha$ on $N_{\infty}$ can be written as 
$$\alpha(t,\sigma)=\beta(t,\sigma)+dt\w\gamma(t,\sigma),$$ 
where, for each fixed $t\in(T,\infty)$, $\beta(t,\sigma)$ and $\gamma(t,\sigma)$ are $m$- and $(m\!-\!1)$-forms on $\Sigma$ respectively.  Therefore, if $\pi:(0,\infty)\times\Sigma\rightarrow\Sigma$ is 
the natural projection, 
 $\Lambda^mT^*N_{\infty}\cong \pi^*(\Lambda^mT^*\Sigma)\oplus\pi^*(\Lambda^{m-1}T^*\Sigma)$. 
Hence, $(d+d^*)_{\infty}$ maps sections of
$\pi^*(\Lambda^2T^*\Sigma)\oplus\pi^*(\Lambda^{\text{odd}}T^*\Sigma)$ to sections of $\pi^*(\Lambda^{\text{odd}}T^*\Sigma)\oplus
\pi^*(\Lambda^{\text{even}}T^*\Sigma)$.  Moreover, this action is given by:
\begin{align}
(d+d^*)_{\infty}&\left(\begin{array}{c}\alpha(t,\sigma) \\ \beta(t,\sigma)+\gamma(t,\sigma)\end{array}\right)=\nonumber\\
&\qquad\left(\begin{array}{cc}\displaystyle d_{\Sigma}+d^*_{\Sigma}&-(\frac{\partial}{\partial t}+3-m) \\
\frac{\partial}{\partial t}+m & -(d_{\Sigma}+d^*_{\Sigma})\end{array}\right)\left(\begin{array}{c}\alpha(t,\sigma) \\ 
\beta(t,\sigma)+\gamma(t,\sigma)\end{array}\right),\label{dd*inftyeq2}
\end{align}
where $m$ denotes the operator which multiplies $m$-forms by a
factor $m$, and $d_{\Sigma}$ and $d^*_{\Sigma}$ are the exterior derivative and its formal adjoint on $\Sigma$.  

However, we wish only to consider elements of
$\Lambda^1T^*\Sigma\oplus\Lambda^2T^*\Sigma$ which correspond, via $\pi^*$, to
self-dual 2-forms on $N_{\infty}$.  Thus we define
$V_{\Sigma}\subseteq\Lambda^2T^*\Sigma\oplus\Lambda^{\text{odd}}T^*\Sigma$
by
\begin{equation}\label{VSigmaeq}
V_{\Sigma}=\{(\alpha,\ast_{\Sigma}\alpha+\beta)\,:\,\alpha\in\Lambda^2T^*\Sigma,\,\beta\in\Lambda^3T^*\Sigma\},
\end{equation}
where $\ast_{\Sigma}$ is the Hodge star on $\Sigma$.  
Then $\pi^*(V_{\Sigma})\cong\Lambda_+^2T^*N_{\infty}\oplus\Lambda^4T^*N_{\infty}$.  We also want to project the
image of $(d+d^*)_{\infty}$ to 3-forms on
$N_{\infty}$, so we let
$$W_{\Sigma}=\{(\beta,\alpha)\,:\,\alpha\in\Lambda^2T^*\Sigma,\, \beta\in\Lambda^3T^*\Sigma\}$$
and let $\pi_{W_{\Sigma}}$ be the projection to $W_{\Sigma}$.  Note that $\pi^*(W_{\Sigma})\cong\Lambda^3T^*N_{\infty}$.

For $w\in\C$, define a map $(d+d^*)_{\infty}(w):C^{1}_{\loc}(V_{\Sigma}\otimes\C)\rightarrow
C^{0}_{\loc}(W_{\Sigma}\otimes\C)$ by:
\begin{align}
(d+d^*)_{\infty}(w)&\left(\begin{array}{c}\alpha(\sigma)\\
*_{\Sigma}\alpha(\sigma)+\beta(\sigma)\end{array}\right)=\nonumber\\
&\pi_{W_{\Sigma}}\circ\left(\begin{array}{cc}\displaystyle d_{\Sigma}+d^*_{\Sigma}&-(w+3-m) \\
w+m & -(d_{\Sigma}+d^*_{\Sigma})\end{array}\right)\left(\begin{array}{c}\alpha(\sigma)\\
*_{\Sigma}\alpha(\sigma)+\beta(\sigma)\end{array}\right).\label{dd*inftyweq1}
\end{align}
 Notice that we have
formally substituted $w$ for $\frac{\partial}{\partial t}$ in
\eq{dd*inftyeq2}.

Let $p>4$ and $k\geq 2$ as in \eq{dd*eq}, noting that $L^p_{k+1}\hookrightarrow C^2$ on $\Sigma$ by the Sobolev Embedding
Theorem. 
Define $\mathcal{C}\subseteq\C$ as the set of $w$ for which the map
\begin{equation}\label{dd*inftyweq2}
(d+d^*)_{\infty}(w):L_{k+1}^p(V_{\Sigma}\otimes\C) \rightarrow L_k^p
(W_{\Sigma}\otimes\C)
\end{equation} 
is \emph{not} an isomorphism. By the proof of \cite[Theorem
1.1]{LockhartMcOwen}, $\mathcal{D}=\{\Re
w:w\in\mathcal{C}\}$. In fact, $\mathcal{C}\subseteq\R$ by
\cite[Lemma 6.1.13]{Marshall}, which shows that the corresponding
sets $\mathcal{C}(\Delta^m)$ are all real.  Hence
$\mathcal{C}=\mathcal{D}$.

The symbol, hence the index $\text{ind}_w$, of $(d+d^*)_{\infty}(w)$
is independent of $w$. Furthermore, $(d+d^*)_{\infty}(w)$ is an
isomorphism for generic values of $w$ since
$\mathcal{D}$ is countable and discrete. Therefore
$\text{ind}_w=0$ for all $w\in\C$; that is,
\begin{equation*}
\dim\Ker(d+d^*)_\infty(w)=\dim\Coker(d+d^*)_\infty(w),
\end{equation*}
so that \eq{dd*inftyweq2} is not an isomorphism precisely when it is not
injective.

The condition $(d+d^*)_{\infty}(w)=0$, using \eq{dd*inftyweq1} and elliptic regularity,
corresponds to the existence of $\alpha\in
C^{\infty}(\Lambda^2T^*\Sigma)$ and $\beta\in
C^{\infty}(\Lambda^3T^*\Sigma)$ satisfying
\begin{align}
\label{Deq}
d_{\Sigma}\alpha=w\beta\quad\text{and}\quad
d_{\Sigma}\!*_{\Sigma}\!\alpha+d^*_{\Sigma}\beta=(w+2)\alpha.
\end{align}
We first note that \eq{Deq} implies that
\begin{equation}\label{betaharmoniceq}
d_{\Sigma}d^*_{\Sigma}\beta=\Delta_{\Sigma}\beta=w(w+2)\beta.
\end{equation}
Since eigenvalues of the
Laplacian on $\Sigma$ are positive, $\beta=0$ if
$w\in(-2,0)$.  If $w=0$ and we take $\alpha=0$, \eq{Deq} forces
$\beta$ to be coclosed. As there are nontrivial coclosed 3-forms on
$\Sigma$, $(d+d^*)_{\infty}(0)$ is not injective and hence
$0\in\mathcal{D}$.
Suppose that $w=-2$ lies in $\mathcal{D}$. Then
\eq{Deq} gives $[\beta]=0$ in $H_{\text{dR}}^3(\Sigma)$. We
know that $\beta$ is harmonic from \eq{betaharmoniceq} so, by Hodge theory, $\beta=0$.
Therefore $-2\in\mathcal{D}$ if and only if there exists
a nonzero closed and coclosed 2-form on $\Sigma$.

\medskip

We state a proposition which follows from the work above.

\begin{prop}\label{Dprop}  Recall the definition of $\Sigma$ from the start of $\S$\ref{defmapsection} and denote 
the Hodge star, the exterior derivative and its formal adjoint on $\Sigma$ by
$*_{\Sigma}$, $d_{\Sigma}$ and $d_{\Sigma}^*$ respectively. 
Let 
$$D(\mu)=\{(\alpha,\,\beta)\in
C^{\infty}(\Lambda^2T^*\Sigma\oplus\Lambda^3T^*\Sigma)
\,:\,d_{\Sigma}\alpha=\mu\beta,\;d_{\Sigma}*_{\Sigma}\alpha+d^*_{\Sigma}\beta=(\mu+2)\alpha\}.$$ 
The countable discrete 
set $\mathcal{D}$ of $\mu\in\R$ for which $(d+d^*)_{\mu}$, given in \eq{dd*eq}, is not Fredholm is given by:
$$\mathcal{D}=\{\mu\in\R\,:\,D({\mu})\neq 0\}.$$
Moreover, $-2\in\mathcal{D}$ if and only if $b^1(\Sigma)>0$, and $0\in\mathcal{D}$.
\end{prop}

\noindent A perhaps more illuminating way to characterise $D(\mu)$
is by:
\begin{align*}
(\alpha,\beta)&\in D(\mu) \Longleftrightarrow\\
&\;\xi=(r^{\mu+2}\alpha+r^{\mu+1}dr\wedge*_{\Sigma}\alpha,
r^{\mu+3}dr\wedge\beta)\in C^{\infty}(\Lambda^2_+T^*C\oplus\Lambda^4T^*C)\\&\;\text{is an $O(r^\mu)$ solution of
$(d+d^*)\xi=0$ in $C^{\infty}(\Lambda^3T^*C)$ as $r\rightarrow\infty$.}
\end{align*}

We now make a definition as in \cite{LockhartMcOwen}.

\begin{dfn}\label{jumpdfn}  Recall the map $(d+d^*)_{\infty}$ defined by \eq{dd*inftyeq2}, the bundle $V_{\Sigma}$
given in \eq{VSigmaeq} and the set $\mathcal{D}$ given in Proposition \ref{Dprop}.  
Let $\mu\in\mathcal{D}$.  We define $\d(\mu)$ to be the dimension
of the vector space of solutions of $(d+d^*)_{\infty}\xi=0$ of
the form
$$\xi(t,\sigma)=e^{\mu t}p\,(t,\sigma),$$
where $p\,(t,\sigma)$ is a polynomial in $t\in(T,\infty)$ with
coefficients in $C^{\infty}(V_{\Sigma}\otimes\C)$. 
\end{dfn}

\noindent The next result is immediate from \cite[Theorem
1.2]{LockhartMcOwen}.
\begin{prop}\label{jumpprop}  Use the notation of Proposition \ref{Dprop} and Definition \ref{jumpdfn}.
Let $\lambda,\lambda^{\prime}\notin\mathcal{D}$ with
$\lambda^{\prime}\leq\lambda$.  For any $\mu\notin\mathcal{D}$ let
$i_{\mu}(d+d^*)$ denote the Fredholm index of the map \eq{dd*eq}.
Then
$$i_{\lambda}(d+d^*)-i_{\lambda^{\prime}}(d+d^*)=\sum_{\mu\in\mathcal{D}\,
\cap(\lambda^{\prime},\,\lambda)}\!\!\!\!\!\!\d(\mu) .$$
\end{prop}

\begin{note} 
A similar result holds for any uniformly elliptic AC operator on $N$.
\end{note}

\noindent We make a key observation, which shall be used on a number of
occasions later.

\begin{prop}\label{nochangeprop}
Use the notation of Theorem \ref{Pexceptthm}.  Let
$\lambda,\lambda^\prime\in\R$ be such that $\lambda^\prime\leq\lambda$
and $[\lambda^\prime,\lambda]\cap\mathcal{D}(\mathcal{P})=\emptyset$.  The
kernels, and cokernels, of $\mathcal{P}:L^p_{k+l,\,\mu}(V)\rightarrow
L^p_{k,\,\mu-\nu}(W)$ when $\mu=\lambda$ and $\mu=\lambda^\prime$
are equal.
\end{prop}

\begin{proof}
Denote the dimensions of the kernel and cokernel of
$\mathcal{P}:L^p_{k+l,\,\mu}(V)\rightarrow L^p_{k,\,\mu-\nu}(W)$, for
$\mu\notin\mathcal{D}(\mathcal{P})$, by $k(\mu)$ and $c(\mu)$ respectively.  Notice that these dimensions 
are finite since $\mathcal{P}$ is Fredholm if $\mu\notin\mathcal{D}(\mathcal{P})$.  
Since $[\lambda^\prime,\lambda]\cap\mathcal{D}(\mathcal{P})=\emptyset$,
$k(\lambda)-c(\lambda)=k(\lambda^\prime)-c(\lambda^\prime)$  by \cite[Theorem
1.2]{LockhartMcOwen} and
hence
\begin{equation}\label{kceq}
 k(\lambda)-k(\lambda^\prime)=c(\lambda)-c(\lambda^\prime).
\end{equation}

We know that $k(\lambda)\geq k(\lambda^\prime)$ because
$L^p_{k+1,\,\lambda^\prime}\hookrightarrow L^p_{k+1,\,\lambda}$ by
Theorem \ref{wSobthm}(a) as $\lambda\geq\lambda^\prime$.
Similarly, since $c(\mu)$ is equal to the dimension of the kernel of
the formal adjoint operator acting on the dual Sobolev space with weight
$-4-(\mu-\nu)$ (as noted in Definition \ref{dualdfn}), $c(\lambda)\leq c(\lambda^\prime)$.  
Since the
left-hand side of \eq{kceq} is non-negative and the
right-hand side is non-positive, we conclude that both
must be zero. The result follows from the fact that
the kernel of $\mathcal{P}$ in $L^p_{k+1,\,\lambda^\prime}$ is contained in
the kernel of $\mathcal{P}$ in $L^p_{k+1,\,\lambda}$, and vice versa for the
cokernels.
\end{proof}

We conclude with an explicit description of the
quantity $\d(\mu)$ for $\mu\in\mathcal{D}$.  This result, as can be seen from the proof, is similar
to \cite[Proposition 2.4]{Joyce5}.

\begin{prop}\label{dmuprop}  In the notation of Proposition \ref{Dprop} and Definition \ref{jumpdfn},
$\d(\mu)=\dim D(\mu)$ for $\mu\in\mathcal{D}$.
\end{prop}

\begin{proof} 
Use the notation of Proposition \ref{Dprop} and the work
preceding it and Definition \ref{jumpdfn}.  Let
$p(t,\sigma)$ be a polynomial in $t\in(T,\infty)$ of degree $m$ written as
$$p(t,\sigma)=\left(\sum_{j=0}^mp_j(\sigma)t^j,\,\sum_{j=0}^m\big(*_{\Sigma}p_j(\sigma)+q_j(\sigma)\big)t^j\right),$$
where $p_j\in C^{\infty}(\Lambda^{2}T^*\Sigma)$ and $q_j\in
C^{\infty}(\Lambda^{3}T^*\Sigma)$ for $j=0,\ldots,m$, with $p_m$ and
$q_m$ not both zero, and let $\xi(t,\sigma)=e^{\mu t}p(t,\sigma)$ as in Definition \ref{jumpdfn}.  We want
to find the dimension $\d(\mu)$ of the space of $\xi$ such that $(d+d^*)_{\infty}\xi=0$.

Using \eq{dd*inftyeq2}, $(d+d^*)_{\infty}\xi=0$ is equivalent to
\begin{align}
\label{pqeq1}
\sum_{j=0}^mt^j(d_{\Sigma}p_j-\mu
q_j)-\sum_{j=0}^mjt^{j-1}q_j&=0\!\quad\!\text{and}
\\[4pt]
\label{pqeq2}
\sum_{j=0}^mt^j\big((\mu+2)p_j-d_{\Sigma}\!*_{\Sigma}\!p_j-d^*_{\Sigma}q_j\big)+\sum_{j=0}^mjt^{j-1}p_j&=0.
\end{align}
Comparing coefficients of $t^m$ we deduce that $(p_m,q_m)\in
D({\mu})$. 

Suppose, for a contradication, that $m\geq 1$.  Comparing coefficients of $t^{m-1}$ in
\eq{pqeq1} and \eq{pqeq2}:
\begin{align}
\label{pqeq3}
d_{\Sigma}p_{m-1}-\mu q_{m-1}=mq_m\quad\text{and}\quad
d_{\Sigma}\!*_{\Sigma}\!p_{m-1}+d^*_{\Sigma}q_{m-1}-(\mu+2)p_{m-1}=mp_m.
\end{align}
We then compute using \eq{pqeq3} and the fact that $(p_m,q_m)\in D(\mu)$:
\begin{align*}
m\langle p_m,p_m\rangle_{L^2}&=\langle
p_m,d_{\Sigma}\!*_{\Sigma}\!p_{m-1}+d^*_{\Sigma}q_{m-1}-(\mu+2)p_{m-1} \rangle_{L^2}\\
&=\langle d_{\Sigma}\!*_{\Sigma}\!p_m-(\mu+2)p_m,p_{m-1}\rangle_{L^2} +\langle
d_{\Sigma}p_m,q_{m-1}\rangle_{L^2} \\
&=\langle -d^*_{\Sigma}q_m,p_{m-1}\rangle_{L^2}+\langle\mu
q_m,q_{m-1}\rangle_{L^2}\\ &=-\langle q_m,d_{\Sigma}p_{m-1}-\mu
q_{m-1}\rangle_{L^2}\\  &=-m\langle q_m,q_m\rangle_{L^2}.
\end{align*}
Hence,
$$m(\|p_m\|_{L^2}^2+\|q_m\|_{L^2}^2)=0$$
and so $p_m=q_m=0$, our required contradiction.  

Thus, the solutions $\xi$ must be of the from
$$\xi(t,\sigma)=e^{\mu t}p(t,\sigma)=e^{\mu t}\big(p_0(\sigma),\,*_{\Sigma}p_0(\sigma)+q_0(\sigma)\big)$$
for $(p_0,q_0)\in D(\mu)$.  
The proposition follows.
\end{proof}

\subsection{The image of {\boldmath $d+d^*$}}



We remarked earlier upon the connection of the kernel of $(d+d^*)_{\mu}$ with the infinitesimal deformations 
of $N$.  This suggests that the \emph{cokernel} of $(d+d^*)_{\mu}$ is 
related to the \emph{obstruction theory} for the deformation problem.  When $(d+d^*)_{\mu}$ is Fredholm, the 
cokernel is isomorphic, via the dual pairing given in Definition \ref{dualdfn}, 
to the kernel of the \emph{adjoint map} which we now define.

\begin{dfn}\label{adjointdfn}
Let $p>4$, $k\geq 2$ and $\mu<1$.  Let $q>1$ be such that $\frac{1}{p}+\frac{1}{q}=1$ and let $l\geq 2$.  
The adjoint map to \eq{dd*eq} is given by
\begin{equation}\label{adjointeq}
(d^*_++d)_{\mu}:L^q_{l+1,\,-3-\mu}(\Lambda^3T^*N)\rightarrow L^q_{l,\,-4-\mu}(\Lambda^2_+T^*N\oplus\Lambda^4T^*N),
\end{equation}
where $(d^*_++d)_{\mu}(\gamma)=d^*_+\gamma+d\gamma$ and $d^*_+=\frac{1}{2}(d^*+*d^*)$.  Let 
$\mathcal{C}_+(\mu)=\Ker (d^*_++d)_{\mu}$.
\end{dfn}

\vspace{-8pt}

\noindent It is straightforward to see, using integration by parts and the dual pairing, 
that $(d^*_++d)_{\mu}$ is the formal adjoint to $(d+d^*)_{\mu}$.

\begin{remark}
The space $\mathcal{C}_+(\mu)$ is the kernel of an elliptic map so, for $\mu\notin\mathcal{D}$, it is a finite-dimensional 
space of smooth forms by Corollary \ref{ACopkerregcor}.  Thus, $\mathcal{C}_+(\mu)$ is independent of $l$ and we can choose $l\geq 2$ 
for use later.
\end{remark}

The aim of this subsection is to prove the following.

\begin{thm}\label{splitthm} Use the notation of Proposition \ref{Dprop} and Definition \ref{adjointdfn}.  Suppose further that
$\mu\in(-\infty,1)\setminus\mathcal{D}$ and let
\begin{equation}\label{fullcokereq}
\mathcal{C}(\mu)=\{\gamma\in L^q_{l+1,\,-3-\mu}(\Lambda^3T^*N)\,:\,d\gamma=d^*\gamma=0\}.
\end{equation}
There exist finite-dimensional subspaces 
$\tilde{\mathcal{C}}(\mu)$ and $\mathcal{O}(N,\mu)$ of $L^p_{k,\,\mu-1}(\Lambda^3T^*N)$ such that
\begin{gather}
L^p_{k,\,\mu-1}(\Lambda^3T^*N)=\Big(\overline{d\big(L^p_{k+1,\,\mu}(\Lambda^2T^*N)\big)}
+ d^*\big(L^p_{k+1,\,\mu}(\Lambda^4T^*N)\big)\Big)\oplus\tilde{\mathcal{C}}(\mu)\label{spliteq1}\\
\intertext{and}
\overline{d\big(L^p_{k+1,\,\mu}(\Lambda^2T^*N)\big)}=
d\big(L^p_{k+1,\,\mu}(\Lambda^2_+T^*N)\big)\oplus\mathcal{O}(N,\mu).\label{spliteq2}
\end{gather}
 Moreover, $\tilde{\mathcal{C}}(\mu)\cong\mathcal{C}(\mu)$ and $\tilde{\mathcal{C}}_+(\mu)=\tilde{\mathcal{C}}(\mu)\oplus\mathcal{O}(N,\mu)\cong\mathcal{C}_+(\mu)$ via the dual pairing.  
\smallskip
Furthermore:
\begin{itemize}
\item[\emph{(a)}] if $\mu<-2$, the sum in \eq{spliteq1} is a  direct sum;
\item[\emph{(b)}] if $\mu\in[-2,0)$, $\mathcal{O}(N,\mu)=0$ and the sum in \eq{spliteq1} is a  direct sum;
\item[\emph{(c)}] if $\mu\in[0,1)$, $\mathcal{O}(N,\mu)=0$ but the sum in \eq{spliteq1} is not necessarily direct.
\end{itemize}
\end{thm}

Before proving the theorem we make an elementary observation.

\begin{lem}\label{maxlem}
If $\mu< 0$, $d\big(L^p_{k+1,\,\mu}(\Lambda^2T^*N)\big)\cap d^*\big(L^p_{k+1,\,\mu}(\Lambda^4T^*N)\big)=\{0\}$.
\end{lem}

\begin{proof}
If $\beta\in L^p_{k+1,\,\mu}(\Lambda^4T^*N)$ such that $d^*\beta$ is exact, then $dd^*\beta=0$.  Therefore $*\beta$ is 
a harmonic function which is $o(\rho^{\mu})$ as $\rho\rightarrow\infty$.  Applying the Maximum Principle, $*\beta=0$.
\end{proof}

\begin{proof}[Theorem \ref{splitthm}(a)] 
The key to proving this part of the theorem lies with comparing the image of $(d+d^*)_{\mu}$ with the image of the map
\begin{equation}\label{dd*eq2}
d+d^*:L^p_{k+1,\,\mu}(\Lambda^2T^*N\oplus\Lambda^4T^*N)\rightarrow L^p_{k,\,\mu-1}(\Lambda^3T^*N)
\end{equation}
given by $(d+d^*)(\alpha,\beta)=d\alpha+d^*\beta$.  Clearly, there exists a finite-dimensional subspace 
$\tilde{\mathcal{C}}(\mu)$ of $L^p_{k,\,\mu-1}(\Lambda^3T^*N)$, which is isomorphic (via the dual pairing) to the 
\emph{annihilator} of the image of \eq{dd*eq2}, such that
\begin{equation}\label{spliteq3}
L^p_{k,\,\mu-1}(\Lambda^3T^*N)=\overline{\Big(d\big(L^p_{k+1,\,\mu}(\Lambda^2T^*N)\big)
\oplus d^*\big(L^p_{k+1,\,\mu}(\Lambda^4T^*N)\big)\Big)}\oplus\tilde{\mathcal{C}}(\mu).
\end{equation}
(The direct sum follows from Lemma \ref{maxlem} as $\mu<-2$.)

Now, the annihilator $\mathcal{A}(\mu)$ of the image of \eq{dd*eq2} is given by:
\begin{align*}
\mathcal{A}(\mu)=\{\gamma\in L^q_{l+1,\,-3-\mu}(&\Lambda^3T^*N)\,:\,\langle \gamma, d\alpha+d^*\beta \rangle=0
\\&\qquad\quad\text{for all $(\alpha,\beta)\in L^p_{k+1,\,\mu}(\Lambda^2T^*N\oplus\Lambda^4T^*N)$}\}.
\end{align*}
Therefore, using integration by parts (justified by the choice of weight for the dual Sobolev space), we deduce that 
$\mathcal{A}(\mu)=\mathcal{C}(\mu)$.

Using the Fredholmness of $(d+d^*)_{\mu}$, since $\mu\notin\mathcal{D}$, there exists a finite-dimensional subspace
 $\tilde{\mathcal{C}}_+(\mu)$ of $L^p_{k,\,\mu-1}(\Lambda^3T^*N)$ such that
\begin{equation}\label{spliteq4}
L^p_{k,\,\mu-1}(\Lambda^3T^*N)=d\big(L^p_{k+1,\,\mu}(\Lambda^2_+T^*N)\big)
\oplus d^*\big(L^p_{k+1,\,\mu}(\Lambda^4T^*N)\big)\oplus\tilde{\mathcal{C}}_+(\mu).
\end{equation}
Moreover, $\tilde{\mathcal{C}}_+(\mu)$ is isomorphic to the annihilator of the image of \eq{dd*eq}, which is equal to 
$\mathcal{C}_+(\mu)$ as given in Definition \ref{adjointdfn} (again using integration by parts).  

Equation \eq{spliteq4} allows
us to deduce that $d^*\big(L^p_{k+1,\,\mu}(\Lambda^4T^*N)\big)$ is closed.  
Hence, Lemma \ref{maxlem} allows us to deduce that $$\overline{d\big(L^p_{k+1,\,\mu}(\Lambda^2T^*N)\big)}\cap
 d^*\big(L^p_{k+1,\,\mu}(\Lambda^4T^*N)\big)=\{0\}.$$  Comparing \eq{spliteq3} and \eq{spliteq4} we see that there exists a 
 finite-dimensional space $\mathcal{O}(N,\mu)$ such that $\tilde{\mathcal{C}}_+(\mu)
 =\tilde{\mathcal{C}}(\mu)\oplus\mathcal{O}(N,\mu)$, from which (a) follows.
\end{proof}

We shall require some preliminary technical results before resuming our proof of Theorem \ref{splitthm}.  We 
consider the maps
\begin{align}\label{d+eq}
&d:L^p_{k+1,\,\mu}(\Lambda^2_+T^*N)\rightarrow L^p_{k,\,\mu-1}(\Lambda^3T^*N)\quad\text{and}\\
\label{deq}
&d:L^p_{k+1,\,\mu}(\Lambda^2T^*N)\rightarrow L^p_{k,\,\mu-1}(\Lambda^3T^*N).
\end{align}

\begin{prop}\label{annprop}  Let $p>4$, $k\geq 2$ and $\mu\in(-\infty,1)\setminus\mathcal{D}$, where $\mathcal{D}$ is given in Proposition 
\ref{Dprop}.  Let the annihilators of the images of \eq{d+eq} and \eq{deq} be $\mathcal{A}_+(\mu)$ and $\mathcal{A}(\mu)$
 respectively.  Then, if $q>1$ such that $\frac{1}{p}+\frac{1}{q}=1$ and $l\geq 2$, 
 \begin{align}
\mathcal{A}_+(\mu)&=
\{\gamma\in L^q_{l+1,\,-3-\mu}(\Lambda^3T^*N)\,:\,d^*\gamma\in L^q_{l,\,-4-\mu}(\Lambda^2_-T^*N)\}\label{A+eq}
\intertext{and}
\mathcal{A}(\mu)&=\{\gamma\in L^q_{l+1,\,-3-\mu}(\Lambda^3T^*N)\,:\,d^*\gamma=0\}.\label{Aeq}
\end{align}
Furthermore, if $\mu\in[-2,1)\setminus\mathcal{D}$, $\mathcal{A}_+(\mu)=\mathcal{A}(\mu)$.
\end{prop}

\begin{proof} The formulae \eq{A+eq} and \eq{Aeq} are easily deduced using the dual pairing and integration by parts.
Clearly $\mathcal{A}(\mu)\subseteq\mathcal{A}_+(\mu)$, so suppose $\gamma\in\mathcal{A}_+(\mu)$ and $\mu\geq -2$.
 Notice that $p>4$ and $\frac{1}{p}+\frac{1}{q}=1$ force $q\in(1,\frac{4}{3})$, 
and $\mu\geq -2$ implies that $-4-\mu\leq -2$.  Therefore, by Theorem \ref{wSobthm}(a), 
$L_{l,\,-4-\mu}^{q}\hookrightarrow L^2_{0,\,-2}=L^2$, recalling \eq{Lpeq}. Thus,
\begin{equation*}
\|d^*\gamma\|_{L^2}^2=\int_N-d^*\gamma\w d^*\gamma =\int_N-d\!*\!\gamma\w
d\!*\!\gamma=\int_N-d(*\gamma\w d\!*\!\gamma)=0.
\end{equation*} 
The integration by parts is valid since $*\gamma=o(\rho^{-3-\mu})$ and $d\!*\!\gamma=o(\rho^{-4-\mu})$ as 
$\rho\rightarrow\infty$, and $-7-2\mu\leq -3$.  The proof is thus complete.
\end{proof}

Using this proposition we can deduce the following invaluable result.

\begin{prop}\label{dimageprop}  Recall the definition of $\mathcal{D}$ in Proposition \ref{Dprop}.
For $p>4$, $k\geq 2$ and $\mu\in [-2,1)\setminus\mathcal{D}$,
$$\overline{d\left(L^p_{k+1,\,\mu}(\Lambda^2T^*N)\right)}=d\left(L^p_{k+1,\,\mu}(\Lambda^2T^*N)\right)=d\left( L^p_{k+1,\,\mu}(\Lambda^2_+T^*N)\right).$$
\end{prop}

\begin{proof}  By Proposition \ref{annprop}, the annihilators of the images of \eq{d+eq} and \eq{deq}
are equal.  We deduce that the \emph{closure} of the ranges of \eq{d+eq} and \eq{deq} are equal.  However, we know that
$\mu\notin\mathcal{D}$ so that \eq{dd*eq} is Fredholm, which means that it has 
closed range.  Therefore, \eq{d+eq} has closed range and hence
$$d\left( L^p_{k+1,\,\lambda}(\Lambda^2_+T^*N)\right)=\overline{d\left(L^p_{k+1,\,\lambda}(\Lambda^2T^*N)\right)}.$$
The left-hand side of this equation is contained in the image of \eq{deq}, whereas the right-hand side contains the image
of \eq{deq}.  The result follows.  
\end{proof}

\begin{proof}[Theorem \ref{splitthm}(b)-(c)]
If $\mu\in[-2,1)\setminus\mathcal{D}$, 
\eq{dd*eq} is Fredholm, so there exists a finite-dimensional subspace  
$\tilde{\mathcal{C}}(\mu)$ of $L^p_{k,\,\mu-1}(\Lambda^3T^*N)$ such that
\begin{equation*}
L^p_{k,\,\mu-1}(\Lambda^3T^*N)=\Big(d\big(L^p_{k+1,\,\mu}(\Lambda^2_+T^*N)\big)+ 
d^*\big(L^p_{k+1,\,\mu}(\Lambda^4T^*N)\big)\Big)\oplus\tilde{\mathcal{C}}(\mu).
\end{equation*}
The results now follow from Lemma \ref{maxlem} and Proposition \ref{dimageprop}.
\end{proof}

\subsection{The kernel of the adjoint map}

We already remarked that studying the kernel of the adjoint map will help us to understand the obstructions to 
deformations of $N$.  For this subsection we use the notation 
of Definition \ref{adjointdfn} and Proposition \ref{annprop}.  Suppose further that $\mu\in[-2,1)\setminus\mathcal{D}$, where $\mathcal{D}$ is given 
in Proposition \ref{Dprop}.

Let $\gamma\in\mathcal{A}_+(\mu)$.  By Proposition \ref{annprop}, $\ast\gamma\in
L_{l+1,\,-3-\mu}^{q}(\Lambda^1T^*N)$ and satisfies
$d\!\ast\!\gamma=0$. Recall that $\Psi:(R,\infty)\times\Sigma\cong N\setminus K$ 
and that $(r,\sigma)$ are the coordinates on $(R,\infty)\times\Sigma$.  So, on $(R,\infty)\times\Sigma$,
$$\ast\gamma\big(\Psi(r,\sigma)\big)=\chi(r,\sigma)+f(r,\sigma)\,dr$$
for a function $f$ and 1-form $\chi$ on $(R,\infty)\times\Sigma$,
where $\chi$ has no $dr$ component. Write the exterior derivative on $(R,\infty)\times\Sigma$ in terms of the exterior derivative
$d_{\Sigma}$ on $\Sigma$ as:
$$d=d_\Sigma+dr\w\frac{\partial}{\partial r}\,.$$
The equation $d\!\ast\!\gamma=0$ then implies that
\begin{equation}\label{d*gammaeq}
d_{\Sigma}\chi=0\qquad\text{and}\qquad
\frac{\partial\chi}{\partial r}-d_{\Sigma}f=0.
\end{equation} 
Define a function $\zeta$ on
$(R,\infty)\times\Sigma$ by
\begin{equation*}
\zeta(r,\sigma)=-\int_r^\infty f(s,\sigma)\,ds.
\end{equation*} This is
well-defined since the modulus of $f$ is $o(r^{-3-\mu})$ as
$r\rightarrow\infty$, where $-3-\mu\leq-1$ since $\mu\geq-2$.
Noting that the modulus of $\chi$ with respect to $g_{\Sigma}$ is
$o(r^{-2-\mu})$ as $r\rightarrow\infty$, with $-2-\mu\leq0$, we
calculate using \eq{d*gammaeq}:
\begin{align*}
d\zeta(r,\sigma)&=-\int_r^\infty d_{\Sigma}f(s,\sigma)\,ds+f(r,\sigma)\,dr=-\int_r^\infty
\frac{\partial\chi}{\partial
r}\,(s,\sigma)\,ds+f(r,\sigma)\,dr \\[4pt]
&=\Big[-\chi(s,\sigma)\Big]_r^{\infty}+f(r,\sigma)\,dr=\chi(r,\sigma)+f(r,\sigma)\,dr=\ast\gamma\big(\Psi(r,\sigma)\big).
\end{align*}
If $\{R\}\times\Sigma$ has a tubular neighbourhood in $N$, which can
be ensured by making $R$ larger if necessary, we can extend $\zeta$
smoothly to a function on $N$. Hence $\zeta\in
L_{l+2,\,-2-\mu}^q(\Lambda^0T^*N)$ with $d\zeta=*\gamma$ on
$N\setminus K$.  This leads us to the following.

\begin{prop}\label{gammahatprop} Use the notation of Definition \ref{adjointdfn} and Propositions \ref{Dprop} and \ref{annprop}.  
Let $\gamma\in\mathcal{A}_+(\mu)$ for $\mu\in[-2,1)\setminus\mathcal{D}$.  There exists $\zeta\in
L_{l+2,\,-2-\mu}^q(\Lambda^0T^*N)$ such that
$*\gamma-d\zeta=\hat{\gamma}$ is a closed compactly supported 1-form. Moreover,
the map $\gamma\mapsto[\hat{\gamma}]$ from $\mathcal{C}_+(\mu)\subseteq\mathcal{A}_+(\mu)$ to
$H^1_{\cs}(N)$ is injective.
\end{prop}

\begin{proof}
Clearly our construction above ensures that $\hat{\gamma}$ is a closed 1-form
which is zero outside of the compact subset $K$ of $N$. Thus
$[\xi]\in H_{\cs}^1(N)$.  Suppose that $\gamma\in\mathcal{C}_+(\mu)$ and $[\hat{\gamma}]=0$. Then
$\xi=d\tilde{\zeta}$ for some function 
$\tilde{\zeta}$ with compact support. 
 Therefore
$$0=d^*\!*\!\gamma=d^*(d\zeta+\hat{\gamma})=d^*d(\zeta+\tilde{\zeta}).$$
Hence $\zeta+\tilde{\zeta}$ is a harmonic function of order $o(\rho^{-2-\mu})$ as $\rho\rightarrow\infty$.  
Since $\mu\geq-2$, the Maximum Principle forces $\zeta+\tilde{\zeta}=0$ and the result follows.
\end{proof}

\noindent It follows from Proposition \ref{dimcokerprop} below that the map from
$\mathcal{C}_+(\mu)$ to $H^1_{\cs}(N)$ is an isomorphism when
$\mu\in[-2,0)\setminus\mathcal{D}$.

\section{The deformation theory}\label{modulispacesection}

In this section we prove our main result, which is a local description of the 
moduli space $\mathcal{M}(N,\lambda)$ of AC coassociative deformations of $N$ with rate $\lambda$.
We remind the reader that we are using notation from the start of $\S$\ref{defmapsection}.  We shall further assume that 
$\lambda\notin\mathcal{D}$, where $\mathcal{D}$ is given by Proposition \ref{Dprop}, and we choose some $p>4$ and $k\geq 2$.

\subsection{Deformations and obstructions}

We begin by identifying the infinitesimal deformation space for our moduli space problem.   

\begin{dfn}\label{infdefspacedfn}
The \emph{infinitesimal deformation space} is
\begin{align*}
\mathcal{I}(N,\lambda)=\{\alpha\in L^p_{k+1,\,\lambda}(\Lambda^2_+T^*N)\,:\,d\alpha=0\}.
\end{align*}
Note that $\mathcal{I}(N,\lambda)$ is a subspace of $\mathcal{K}(\lambda)$, given in Definition \ref{dd*dfn}, and is 
thus finite-dimensional as $\lambda\notin\mathcal{D}$.  In fact, $\mathcal{I}(N,\lambda)\cong\mathcal{K}(\lambda)$ if
$\lambda<0$ by Lemma \ref{maxlem}.

We say that our deformation theory is \emph{unobstructed} if $\mathcal{M}(N,\lambda)$, given in Definition \ref{modulispacedfn}, 
is a smooth manifold near $N$ of dimension $\dim\mathcal{I}(N,\lambda)$.
\end{dfn}

By Proposition \ref{Fmodulispaceprop} and Corollary \ref{Fregcor}, $\mathcal{M}(N,\lambda)$ is homeomorphic near $N$ to 
the kernel of the map $F$, given in Definition \ref{Fdfn}, near $0$ in $L^p_{k+1,\,\lambda}(U)$.  Clearly
 $\mathcal{I}(N,\lambda)$ is the tangent space to $\Ker F$ at $0$, so it can be identified with the infinitesimal deformation 
 space.  

Furthermore, if there are no obstructions to the deformation theory of $N$, then every infinitesimal deformation should extend to 
a genuine deformation and the moduli space should be a smooth manifold.   This justifies our definition of an unobstructed 
deformation theory.

\smallskip

Theorem \ref{splitthm} identifies the obstruction space.

\begin{dfn}\label{obsspacedfn} 
The \emph{obstruction space} is the finite-dimensional space 
$$\mathcal{O}(N,\lambda)\cong
 \frac{\overline{d\big(L^p_{k+1,\,\lambda}(\Lambda^2T^*N)\big)}}{d\big(L^p_{k+1,\,\lambda}(\Lambda^2_+T^*N)\big)}$$
given in \eq{spliteq2}.  We observe that
$$\dim\mathcal{O}(N,\lambda)=\dim\mathcal{C}_+(\lambda)-\dim\mathcal{C}(\lambda),$$
in the notation of Definition \ref{adjointdfn} and Theorem \ref{splitthm}. 
\end{dfn}

\noindent It shall become clear in the next subsection why we should think of $\mathcal{O}(N,\lambda)$ as the 
obstructions to our deformation theory.  However, we notice by Theorem \ref{splitthm} that $\mathcal{O}(N,\lambda)=\{0\}$ for 
$\lambda\in[-2,1)$, which should mean that our deformation theory is unobstructed for these rates -- we shall confirm that 
this is the case.

\smallskip

The key step in understanding the obstruction theory is contained in the next result, which studies the image of the
deformation map $F$.

\begin{prop}\label{Fimageprop}
Use the notation of Proposition \ref{tubenbdprop} and Definition \ref{Fdfn}. 
Making the open set $U$ smaller if necessary,
$$F:L^p_{k+1,\,\lambda}(U)\rightarrow \overline{d\big(L^p_{k+1,\,\lambda}(\Lambda^2T^*N)\big)}\subseteq L^p_{k,\,\lambda-1}
(\Lambda^3T^*N).$$
\end{prop}

\begin{proof}
It was noted in the proof of Proposition \ref{isokernelsprop} that $F(\alpha)$ is exact for $\alpha\in L^p_{k+1,\,\lambda}(U)$.  
However, we need to know that we can choose a 2-form $H(\alpha)$ lying in an appropriate weighted Sobolev space such that
 $d\big(H(\alpha)\big)=F(\alpha)$.

Let $u$ be the vector field given by dilations, which, in
coordinates $(x_1,\ldots,x_7)$ on $\R^7$, is written:
\begin{equation}\label{udilationeq}
u=x_1\frac{\partial}{\partial
x_1}+\ldots+x_7\frac{\partial}{\partial x_7}\,.
\end{equation} 
Then
the Lie derivative of $\varphi$ along $u$ is:
\begin{equation}\label{uLiederiveq}
\mathcal{L}_u\varphi=d(u\cdot\varphi)=3\varphi.
\end{equation} 
Therefore,
$\psi=\frac{1}{3}\,u \cdot\varphi$ is a 2-form such that
$d\psi=\varphi$. Note that $\psi|_C\equiv 0$ since
$$(u\cdot\varphi)|_C=u\cdot(\varphi|_C)=0,$$
as $u\in TC$ and $C$ is coassociative by Proposition \ref{coassconeprop}. 
Recall the definition of $f_{\alpha}$ and $N_{\alpha}$ in Definition \ref{Fdfn}.
Define, for $\alpha\in
C^{1}_{\loc}(U)$, 
$$H(\alpha)=f_\alpha^*\big(\psi|_{N_\alpha}\big)$$ 
so that
$$F(\alpha)=f_\alpha^*\big(d\psi|_{N_\alpha}\big)=d\Big(f_\alpha^*\big(\psi|_{N_\alpha}\big)\Big)=d\big(H(\alpha)\big).$$

Recall the diffeomorphism
$\Psi_\alpha:(R,\infty)\times\Sigma\rightarrow N_\alpha\setminus
K_\alpha$, where $K_\alpha$ is compact, introduced before
Proposition \ref{Fmodulispaceprop}, and the inclusion map $\iota:(R,\infty)\times\Sigma\rightarrow C$ 
given by $\iota(r,\sigma)=r\sigma$.  The decay of $H(\alpha)$ at
infinity is determined by:
$$\Psi_\alpha^*(\psi)=(\Psi_\alpha^*-\iota^*)(\psi)+\iota^*(\psi)=(\Psi_\alpha^*-\iota^*)(\psi)$$
since $\psi|_C\equiv0$.  For $(r,\sigma)\in (R,\infty)\times\Sigma$,
\begin{align}
(\Psi_\alpha^*-\iota^*)(\psi)|_{(r,\sigma)}=&
\left(d\Psi_\alpha|^*_{(r,\sigma)}\big(\psi|_{\Psi_\alpha(r,\sigma)}\big)
-d\iota|^*_{(r,\sigma)}\big(\psi|_{\Psi_\alpha(r,\sigma)}\big)\right)\nonumber\\[4pt]
\label{psiesteq} &+
d\iota|^*_{(r,\sigma)}\big(\psi|_{\Psi_\alpha(r,\sigma)}-\psi|_{r\sigma}\big),
\end{align}
using the linearity of $d\iota^*$ to derive the last term.  Since
$|\psi|=O(r)$ and $\Psi_\alpha$ satisfies \eq{Psieq1} so that
$|d\Psi_\alpha^*-d\iota^*|=O(r^{\lambda-1})$ as
$r\rightarrow\infty$, the expression in brackets in \eq{psiesteq} is
$O(r^\lambda)$. The magnitude of the final term in \eq{psiesteq} at infinity is determined by the
behaviour of $d\iota^*$, $\nabla\psi$ and $\Psi_\alpha-\iota$.
Hence, as $|d\iota^*|$ and $|\nabla\psi|$ are $O(1)$, using
\eq{Psieq1} again implies that this term is $O(r^\lambda)$. We
conclude that if $\alpha\in L^p_{k+1,\,\lambda}(U)$ then
$H(\alpha)\in L^p_{k,\,\lambda}(\Lambda^2T^*N)$. Notice that
$H(\alpha)$ has one degree of differentiability less than one would
expect since it depends on $\alpha$ and $\nabla\alpha$.

By Proposition \ref{annprop}, the annihilator $\mathcal{A}(\lambda)$ of the image of \eq{deq} for $\mu=\lambda$ comprises 
of coclosed forms and lies in $L^q_{l+1,\,-3-\lambda}$, where $\frac{1}{p}+\frac{1}{q}=1$ and $l\geq 2$.  
For $\gamma\in\mathcal{A}(\lambda)$, recalling the dual pairing given in Definition \ref{dualdfn}, 
$$\langle F(\alpha),\gamma\rangle=\langle d\big(H(\alpha)\big),\gamma \rangle=
\langle H(\alpha),d^*\gamma \rangle=0,$$
where the integration by parts is valid as $H(\alpha)\in L^p_{k,\,\lambda}$ and $\gamma\in L^q_{l+1,\,-3-\lambda}$.  Thus,
$F(\alpha)$ must lie in the closure of the image of \eq{deq}.  
\end{proof}

\subsection{The moduli space {\boldmath $\mathcal{M}(N,\lambda)$}}

We now state and prove our main result.  

\begin{thm}\label{modulispacethm}  
Let $N$ be a coassociative 4-fold in $\R^7$ which is asymptotically conical 
with rate $\lambda$.  Let 
$p>4$, $k\geq 2$ and suppose that $\lambda\notin\mathcal{D}$, where $\mathcal{D}$ is given in Proposition \ref{Dprop}.  
Use the notation of Definitions \ref{modulispacedfn}, \ref{dd*dfn}, \ref{infdefspacedfn} and \ref{obsspacedfn} and let
$$\mathcal{B}(\lambda)=\big\{\beta\in L^{p}_{k+1,\,\lambda}(\Lambda^4T^*N)
:d^*\beta\in \overline{d\big(L^{p}_{k+1,\,\lambda}(\Lambda^2T^*N)\big)}\subseteq L^p_{k,\,\lambda-1}(\Lambda^3T^*N)\big\}.$$   

  There exist a manifold $\hat{\mathcal{M}}(N,\lambda)$, which
is an open neighbourhood of $0$ in $\mathcal{K}(\lambda)$, and a smooth map 
$\pi:\hat{\mathcal{M}}(N,\lambda)\rightarrow\mathcal{O}(N,\lambda)$, with $\pi(0)=0$, such that an open neighbourhood of
$0$ in $\Ker\pi$ is homeomorphic to 
$\mathcal{M}(N,\lambda)$ near $N$.  

Moreover, if $\mathcal{O}(N,\lambda)=\{0\}$, $\mathcal{M}(N,\lambda)$ is a smooth manifold near $N$ with
$$\dim\mathcal{M}(N,\lambda)=\dim\mathcal{I}(N,\lambda)=\dim\mathcal{K}(\lambda)-\dim\mathcal{B}(\lambda);$$
that is, the deformation theory is unobstructed.
\end{thm}

\begin{note}
By Theorem \ref{splitthm}, $\mathcal{B}(\lambda)=\{0\}$ if $\lambda<0$ and, if $\lambda\geq 0$, 
$$\mathcal{B}(\lambda)=\big\{\beta\in L^{p}_{k+1,\,\lambda}(\Lambda^4T^*N)
:d^*\beta\in d\big(L^{p}_{k+1,\,\lambda}(\Lambda^2_+T^*N)\big)\big\}.$$  
Moreover, $\mathcal{B}(\lambda)$ is finite-dimensional as it is isomorphic to a subspace of the harmonic functions in $L^p_{k+1,\,\lambda}$.
\end{note}

\begin{proof}
Recall the open set $U$ given by Propostion \ref{tubenbdprop}.  Make $U$ smaller if necessary so that Proposition 
\ref{Fimageprop} holds and that the $C^1_1$-norm of elements $\alpha\in L^p_{k+1,\,\lambda}(U)$ is small enough for Corollary 
\ref{Fregcor} to apply.  Define
\begin{align*}
V&=L^p_{k+1,\,\lambda}(U)\times L^p_{k+1,\,\lambda}(\Lambda^4T^*N)\\
X&=L^p_{k+1,\,\lambda}(\Lambda^2_+T^*N)\oplus L^p_{k+1,\,\lambda}(\Lambda^4T^*N)\\
Y&=\mathcal{O}(N,\lambda)\subseteq L^p_{k,\,\lambda-1}(\Lambda^3T^*N)\quad\text{and}\\
Z&=\overline{d\big(L^p_{k+1,\,\lambda}(\Lambda^2T^*N)\big)}+ d^*\big(L^p_{k+1,\,\lambda}(\Lambda^4T^*N)\big)
\subseteq L^p_{k,\,\lambda-1}(\Lambda^3T^*N).
\end{align*}
By Definition \ref{wSobdfn}, $X$ is a Banach space and, as noted at the start of $\S$\ref{FGsubsection}, $V$ is an open neighbourhood
 of zero  in $X$.  The obstruction space $Y$ is trivially a Banach space since it is finite-dimensional. 
Finally, $Z$ is a Banach space as it is a closed subset of a Banach space by the proof of Theorem \ref{splitthm}.  Define a smooth map 
$\mathcal{G}$ on the open subset $V\times Y$ of $X\times Y$ by 
$$\mathcal{G}(\alpha,\beta,\gamma)=G(\alpha,\beta)+\gamma,$$
where $G$ is given in Definition \ref{Gdfn}.  By Proposition \ref{Fimageprop}, 
$\mathcal{G}$ maps into $Z$.  
  Moreover, the linearisation of $\mathcal{G}$ at zero, $d\mathcal{G}|_{(0,0,0)}:X\times Y\rightarrow Z$,
acts as
$$(\alpha,\beta,\gamma)\longmapsto d\alpha+d^*\beta+\gamma.$$
Therefore, $d\mathcal{G}|_{(0,0,0)}$ is surjective by Theorem \ref{splitthm}.  Furthermore, using the notation of
Definition \ref{dd*dfn}, we see that
\begin{align*}
\Ker d\mathcal{G}|_{(0,0,0)}&=\{(\alpha,\beta,\gamma)\in X\times Y\,:\,d\alpha+d^*\beta+\gamma=0\}\\
&\cong\Ker (d+d^*)_{\lambda}=\mathcal{K}(\lambda)
\end{align*}
since, by definition of $Y$, $(d+d^*)_{\lambda}(X)\cap Y=\{0\}$.
Finally note that there exists a closed subspace $A$ of $X$ such that
$X=\mathcal{K}(\lambda)\oplus A$ as $\mathcal{K}(\lambda)$ is the kernel of a Fredholm operator.  

We now apply the Implicit Function Theorem (Theorem \ref{IFthm}) to $\mathcal{G}$.  Thus, there exist open sets
$\hat{\mathcal{M}}(N,\lambda)\subseteq \mathcal{K}(\lambda)$, $W_A\subseteq A$ and $W_Y\subseteq Y$, each containing
zero, with $\hat{\mathcal{M}}(N,\lambda)\times W_A\subseteq V$, and smooth maps $\mathcal{W}_A:\hat{\mathcal{M}}(N,\lambda)
\rightarrow W_A$ and $\mathcal{W}_Y:\hat{\mathcal{M}}(N,\lambda)\rightarrow W_Y$ such that
\begin{align*}
\Ker\mathcal{G}\cap\big(\hat{\mathcal{M}}(N,\lambda)&\times W_A\times W_Y\big)\\
&
=\big\{\big(\xi,\mathcal{W}_A(\xi),\mathcal{W}_Y(\xi)\big)\in \mathcal{K}(\lambda)\oplus 
A\oplus Y\,:\,\xi\in\hat{\mathcal{M}}(N,\lambda)\big\}.
\end{align*}
We take our map $\pi=\mathcal{W}_Y$.  Thus, an open neighbourhood of zero in $\Ker\pi$ is homeomorphic to an open
neighbourhood of zero in $\Ker G\subseteq V$.  

If $\lambda<0$, Proposition \ref{isokernelsprop} and Theorem \ref{splitthm}(a)-(b) imply that $\Ker G\subseteq V$ is isomorphic to 
$\Ker F\subseteq L^p_{k+1,\,\lambda}(U)$, $\mathcal{I}(N,\lambda)\cong\mathcal{K}(\lambda)$ and $\mathcal{B}(\lambda)=\{0\}$.  
Since $\Ker F$ consists of smooth forms by Corollary \ref{Fregcor} and gives a local description of the moduli space 
by Proposition \ref{Fmodulispaceprop}, the result follows for $\lambda<0$.

Therefore, suppose $\lambda\geq 0$.  Define a smooth map $\pi_G$ on $\Ker G$ by $\pi_G(\alpha,\beta)=\beta$.  Notice that
$$\Ker\pi_G\cong\Ker F=\{\alpha\in L^p_{k+1,\,\lambda}(U)\,:\,F(\alpha)=0\}.$$
Moreover, if $(\alpha,\beta)\in\Ker G$, 
$$d^*\beta=-F(\alpha)\in \overline{d\big(L^p_{k+1,\,\lambda}(\Lambda^2T^*N)\big)}$$
by Proposition \ref{Fimageprop}.  Therefore, $\pi_G:\Ker G\rightarrow\mathcal{B}(\lambda)$ and
$d\pi_G|_{(0,0)}:\mathcal{K}(\lambda)\rightarrow \mathcal{B}(\lambda)$.  Since $d\pi_G|_{(0,0)}$ is clearly surjective, $\pi_G$ is locally surjective.  Moreover, 
$$\Ker d\pi_G|_{(0,0)}\cong\{\alpha\in L^p_{k+1,\,\lambda}(\Lambda^2_+T^*N)\,:\,d\alpha=0\}=\mathcal{I}(N,\lambda)=\Ker dF|_{0}.$$ 
Hence, 
$\Ker F$ is locally diffeomorphic to $\Ker d\pi_G|_{(0,0)}$ and the result follows. 
\end{proof}

Theorem \ref{splitthm} gives us an immediate corollary to Theorem \ref{modulispacethm}.

\begin{cor}\label{modulispacecor}
The deformation theory of AC coassociative 4-folds $N$ in $\R^7$ with generic rate $\lambda\in[-2,1)$ is unobstructed.
\end{cor}

\begin{remark}
We have no guarantee that the deformation theory is unobstructed for generic rates $\lambda<-2$.  However, we would expect 
the obstruction space to be zero for generic choices of $N$.
\end{remark}

\section{The dimension of the moduli space {\boldmath $\mathcal{M}(N,\lambda)$}}\label{dimsection}

We shall continue to use the notation from the start of $\S$\ref{defmapsection}.  
We additionally recall the notation of Definitions \ref{dd*dfn} and \ref{adjointdfn}:
for $\mu\in\R$, we denote the kernel of $(d+d^*)_{\mu}$, given in \eq{dd*eq}, by $\mathcal{K}(\mu)$ and 
the kernel of $(d^*_++d)_{\mu}$, given in 
\eq{adjointeq}, by $\mathcal{C}_+(\mu)$.  
We also remind the reader that $\mathcal{C}_+(\mu)$ is isomorphic to the cokernel of \eq{dd*eq}  
via the pairing given in Definition \ref{dualdfn} when $\mu\notin\mathcal{D}$, where $\mathcal{D}$ is given in 
Proposition \ref{Dprop}.  We shall split our study of the dimension of the moduli space into three ranges of rates 
$\lambda\notin\mathcal{D}$: $\lambda\in[-2,0)$, $\lambda\in[0,1)$ and $\lambda<-2$.

\subsection{Topological calculations}

The first proposition we state follows from standard results in
algebraic topology if we consider our coassociative 4-fold $N$, which is AC to $C\cong (0,\infty)\times\Sigma$, 
as the interior of a manifold which has boundary $\Sigma$.  This result is invaluable for our dimension calculations.

\begin{prop}\label{seqprop}
Recall the notation from the start of $\S$\ref{defmapsection}.  
Let the map $\phi_m:H_{\cs}^m(N)\rightarrow
H_{\dR}^m(N)$ be defined by $\phi_m([\xi])=[\xi]$. Let
$r>R$ and let $\Psi_r:\Sigma\rightarrow N$ be the embedding given by
$\Psi_r(\sigma)=\Psi(r,\sigma)$.  Define
$p_m:H_{\dR}^m(N)\rightarrow
H_{\dR}^m(\Sigma)$ by $p_m([\xi])=[\Psi_r^*\xi]$.  Let
$f\in C^{\infty}(N)$ be such that $f=0$ on $K$ and $f=1$ on
$(R+1,\infty)\times\Sigma$.  If
$\pi_\Sigma:(R,\infty)\times\Sigma\cong N\setminus K
\rightarrow\Sigma$ is the projection map, define
$\partial_m:H_{\dR}^m(\Sigma)\rightarrow
H_{\cs}^{m+1}(N)$ by
$\partial_m([\xi])=[d(f\pi_\Sigma^*\xi)]$. Then the following
sequence is exact:
\begin{equation}\label{seqeq}
\cdots\longrightarrow H^m_{\cs}(N)
\,{\buildrel\phi_m\over\longrightarrow}\,
H^m_{\dR}(N)\,{\buildrel p_m\over\longrightarrow}\,
H^m_{\dR}(\Sigma)\,{\buildrel
\partial_m\over\longrightarrow}\,
H_{\cs}^{m+1}(N)\longrightarrow\cdots.\end{equation}
\end{prop}

\begin{remark} The fact that $H^0_{\cs}(N)=H^4_{\dR}(N)=\{0\}$
enables us to calculate the dimension of various spaces more
easily using the long exact sequence \eq{seqeq}.
\end{remark}

We now identify a space of forms which we shall relate to the topology of $N$.

\begin{dfn}\label{Hmdfn}
Define $\mathcal{H}^m(N)$ by 
\begin{equation}\label{Hmeq}
\mathcal{H}^m(N)=\{\xi\in L^2(\Lambda^mT^*N)\,:\,d\xi=d^*\xi=0\}.
\end{equation}
Using \eq{Lpeq}, we notice that $L^2_{0,\,-2}=L^2$.  Moreover, the elliptic regularity 
result Theorem \ref{ACopregthm}(a) applied to the Laplacian can be improved as in \cite[Proposition 5.3]{Lockhart}, using more
theory of weighted Sobolev spaces than we wish to discuss here, to show that harmonic forms in $L^2_{0,\,-2}$, hence
elements of $\mathcal{H}^m(N)$, lie in $L^2_{k,\,-2}$ for all $k\in\N$ and thus are smooth.  
\end{dfn}

The moduli space of deformations for a \emph{compact} coassociative 4-fold $P$, by \cite[Theorem 4.5]{McLean}, is smooth of dimension 
$b^2_+(P)$.  This suggests that we need to define an analogue of $b^2_+$ for an AC coassociative 4-fold $N$.

\begin{dfn}\label{Hplusdfn}
As noted in Definition \ref{Hmdfn}, $\mathcal{H}^2(N)$ consists of smooth forms.  
The Hodge star maps $\mathcal{H}^2(N)$ into itself, so there is a
splitting $\mathcal{H}^2(N)=\mathcal{H}_+^2(N)\oplus\mathcal{H}_-^2(N)$ where
$$\mathcal{H}_{\pm}^2(N)=\mathcal{H}^2(N)\cap
C^\infty(\Lambda_{\pm}^2T^*N).$$

Define $\mathcal{J}(N)=\phi_2\left(H_{\cs}^2(N)\right)$, with $\phi_2$ as in Proposition \ref{seqprop}.  If
$[\alpha]$, $[\beta]\in \mathcal{J}(N)$, there exist compactly supported closed
2-forms $\xi$ and $\eta$ such that $[\alpha]=\phi_2([\xi])$ and $[\beta]=\phi_2([\eta])$.
We define a product on $\mathcal{J}(N)\times \mathcal{J}(N)$ by
\begin{equation}\label{producteq}
[\alpha]\cup[\beta]=\int_N\xi\w\eta. \end{equation} Suppose that
$\xi^{\prime}$ and $\eta^{\prime}$ are also compactly supported with
$[\alpha]=\phi_2([\xi^{\prime}])$ and $[\beta]=\phi_2([\eta^{\prime}])$.  Then there
exist 1-forms $\chi$ and $\zeta$ such that $\xi-\xi^{\prime}=d\chi$ and
$\eta-\eta^{\prime}=d\zeta$. Therefore,
\begin{align*}
\int_N\xi^{\prime}\w\eta^{\prime}&=\int_N(\xi-d\chi)\w(\eta-d\zeta)
=\int_N\xi\w\eta-d\chi\w\eta-\xi^{\prime}\w d\zeta\\[4pt]
&=\int_N\xi\w\eta-d(\chi\w\eta)-d(\xi^{\prime}\w\zeta)=\int_N\xi\w\eta,
\end{align*}
as both $\chi\w\eta$ and $\xi^{\prime}\w\zeta$ have compact support.
The product \eq{producteq} on $\mathcal{J}(N)\times \mathcal{J}(N)$ is thus
well-defined and is a symmetric topological product with a signature
$(a,b)$.  By \cite[Example 0.15]{Lockhart}, $\mathcal{H}^2(N)\cong\mathcal{J}(N)$ via $\alpha\mapsto[\alpha]$.  We therefore define 
\begin{equation*}\label{b2plus}
b^2_+(N)=\dim\mathcal{H}^2_+(N).
\end{equation*}
Hence, $b_+^2(N)=a$, which is a topological number.
\end{dfn}

We can now prove a useful fact about the kernel of $(d+d^*)_{\mu}$.

\begin{lem}\label{dimkerlem} In the notation of Definitions \ref{dd*dfn} and \ref{Hplusdfn}, 
$$\dim\mathcal{K}(\mu)\leq b^2_+(N)\leq\dim\mathcal{K}(-2)$$
for all $\mu<-2$ such that $\mu\notin\mathcal{D}$.
\end{lem}

\begin{proof}
By Theorem \ref{wSobthm}(a), recalling that $p>4$ and $k\geq 2$, $L^p_{k+1,\,\mu}\hookrightarrow L^2_{0,\,-2}$ and
$L^2_{k+3,\,-2}\hookrightarrow L^p_{k+1,\,-2}$.  The result follows from observations in Definitions \ref{infdefspacedfn} and \ref{Hplusdfn}.
\end{proof}

We want to make a similar statement about the cokernel of $(d+d^*)_{\mu}$.  In fact, we can do a lot better, but
we first need a result concerning functions on cones \cite[Lemma 2.3]{Joyce5}.

\begin{prop}\label{conefunctionprop}
Recall the cone $C$ and its link $\Sigma$.  
Suppose that $f:C\rightarrow\R$ is a nonzero function such that
$$f(r,\sigma)=r^{\mu}f_{\Sigma}(\sigma)$$
for some function $f_{\Sigma}:\Sigma\rightarrow\R$ and $\mu\in\R$. Denoting
the Laplacians on $C$ and $\Sigma$ by $\Delta_C$ and
$\Delta_{\Sigma}$ respectively,
$$\Delta_C
f(r,\sigma)=r^{\mu-2}\big(\Delta_{\Sigma}f_{\Sigma}(\sigma)-\mu(\mu+2)f_{\Sigma}(\sigma)\big).$$
Therefore, since $\Sigma$ is compact, $\mu(\mu+2)\geq 0$ and so
there exist no nonzero homogeneous harmonic functions of order
$O(r^{\mu})$ on $C$ with $\mu\in (-2,0)$.  
\end{prop}

\begin{prop}\label{dimcokerprop}  In the notation of 
Definition \ref{adjointdfn} and Theorem \ref{splitthm},
$$\dim\mathcal{C}_+(\mu)=\dim\mathcal{C}(\mu)=b^3(N)$$ 
for all $\mu\in(-2,0)\setminus\mathcal{D}$, and so they are independent of $\mu$ in this range.
\end{prop}

\begin{proof} 
Let $\mu\in(-2,0)\setminus\mathcal{D}$ and let $\gamma\in\mathcal{C}_+(\mu)\subseteq L^q_{l+1,\,-3-\mu}(\Lambda^3T^*N)$. 
 Recall the closed compactly supported 
1-form $\hat{\gamma}=*\gamma-d\zeta$ related to $\gamma$ as given by Proposition \ref{gammahatprop}.  
Clearly $\Delta\zeta=-d^*\hat{\gamma}$ since $d^*\!*\!\gamma=0$. 
Hence, by Proposition \ref{gammahatprop}, $d^*\hat{\gamma}$
lies in the image of the map
$$\Delta_{-2-\mu}=\Delta:L_{l+2,\,-2-\mu}^q(\Lambda^0T^*N)\rightarrow
L_{l,\,-4-\mu}^q(\Lambda^0T^*N).$$  Therefore $d^*\hat{\gamma}$ is
orthogonal, via the pairing given in Definition \ref{dualdfn}, to the kernel of the adjoint of
$\Delta_{-2-\mu}$. Let $\nu\in(-2,0)\setminus\mathcal{D}$.
 Proposition \ref{conefunctionprop} implies, as shown in \cite[Definition 2.5 \& Theorem 2.11]{Joyce5}, that there
are no elements of $\mathcal{D}(\Delta)$, defined in Theorem
\ref{Pexceptthm}, between $-2-\mu$ and $-2-\nu$.  Thus
$\Coker\Delta_{-2-\mu}=\Coker\Delta_{-2-\nu}$.
Obviously, $d^*\hat{\gamma}$ is then orthogonal to the kernel of the
adjoint of $\Delta_{-2-\nu}$.  Consequently, there exists
$\zeta_{\nu}\in L_{l+2,\,-2-\nu}^q(\Lambda^0T^*N)$ such that
$\Delta_{-2-\nu}\zeta_{\nu}=-d^*\hat{\gamma}$ and hence $\zeta-\zeta_\nu$ is
harmonic. Moreover,
$\zeta-\zeta_\mu=o(\rho^{-\min(\mu\,,\,\nu)-2})$ as
$\rho\rightarrow\infty$, where $\rho$ is a radius function on $N$. Since $-\min(\mu,\nu)-2<0$,
we use the Maximum Principle to deduce that $\zeta-\zeta_\mu=0$.
Hence $*\gamma$ and $\gamma$ lie in $L_{l+1,\,-\nu-3}^q$ for any
$\nu\in(-2,0)$. The dimension of $\mathcal{C}_+(\mu)$ is
therefore constant for $\mu\in(-2,0)\setminus\mathcal{D}$. 

We now note that, by \cite[Example 0.15]{Lockhart}, 
$\mathcal{H}^3(N)\cong H^3_{\dR}(N)$ via $\eta\mapsto [\eta]$, where $\mathcal{H}^3(N)$ is given by
\eq{Hmeq}.  By Theorem \ref{wSobthm}(a), 
 $L^q_{l+1,\,-3-\mu}\hookrightarrow L^2_{0,\,-2}=L^2$ for $\mu\geq -1$ and $L^2_{l+1,\,-2}\hookrightarrow L^q_{l+1,\,-3-\mu}$ if 
 $\mu<-1$.  By Theorem \ref{splitthm}, $\mathcal{C}_+(\mu)=\mathcal{C}(\mu)$. 
Hence, $\mathcal{C}_+(\mu)\subseteq\mathcal{H}^3(N)$ whenever $\mu\in [-1,0)\setminus\mathcal{D}$,
and $\mathcal{H}^3(N)\subseteq\mathcal{C}_+(\mu)$ if $\mu\in (-2,-1)\setminus\mathcal{D}$.  The result follows. 
\end{proof}

\subsection{Rates {\boldmath $\lambda\in [-2,0)\setminus\mathcal{D}$}}

For convenience we introduce some extra notation.

\begin{dfn}\label{l-l+dfn}
Since $\mathcal{D}$ is discrete by Proposition \ref{Dprop}, we may choose
$\lambda_-<-2$ and $\lambda_+\in(-2,0)$ such that $[\lambda_-,\lambda_+]\cap\mathcal{D}\subseteq\{-2\}$.
\end{dfn}

\noindent We are thus able to make the following
proposition, which determines the dimension of $\mathcal{M}(N,\lambda)$ in a special case.

\begin{prop}\label{dimprop0}  
Use the notation from Definitions \ref{jumpdfn} and \ref{Hplusdfn}.
If $b^1(\Sigma)=0$ and $\lambda\in [-2,0)\setminus\mathcal{D}$, 
$$\dim\mathcal{M}(N,\lambda)=b^2_+(N)+\!\!\!\!\!\!\sum_{\mu\in\mathcal{D}\,
\cap(-2,\,\lambda)}\!\!\!\!\!\!\d(\mu).$$
\end{prop}

\begin{proof}  Use the notation of Definition \ref{l-l+dfn} and recall that $-2\in\mathcal{D}$ if and only if
$b^1(\Sigma)> 0$ by Proposition \ref{Dprop}.  Since $-2\notin\mathcal{D}$,
$[\lambda_-,\lambda_+]\cap\mathcal{D}=\emptyset$.  By Proposition \ref{nochangeprop}, the dimensions of the kernel
and cokernel of \eq{dd*eq} are constant for $\mu\in[\lambda_-,\lambda_+]$.  Therefore, by
Proposition \ref{dimkerlem}, 
$$\dim\mathcal{K}(\mu)=
\dim\mathcal{K}(-2)=b^2_+(N) \qquad\text{for
$\mu\in[\lambda_-,\lambda_+]$.}$$ 
Applying Propositions \ref{jumpprop} and \ref{dimcokerprop} completes the proof.
\end{proof}

The next proposition enables us to calculate the dimension of
$\mathcal{M}(N,{\lambda})$ when $-2\in\mathcal{D}$.  The proof is long, but the main idea is to identify the forms which 
are added to the kernel of \eq{dd*eq} as the rate crosses $-2$.  By the work in \cite{LockhartMcOwen},
these forms are precisely those which are asymptotic to homogeneous solutions $\xi\in C^{\infty}
(\Lambda^2_+T^*C\oplus\Lambda^4T^*C)$ of order $O(r^{-2})$ to $(d+d^*)\xi=0$, as described after Proposition \ref{Dprop}.   
In the course of the proof, we will need an elementary lemma.

\begin{lem}\label{dimagelem} If $P$ is a compact 4-dimensional Riemannian manifold,  
$$d\big(C^{\infty}(\Lambda^2T^*P)\big)=d\big(C^{\infty}(\Lambda^2_+T^*P)\big).$$
\end{lem}

\begin{proof}
Let $\theta\in C^{\infty}(\Lambda^2T^*P)$.  By Hodge theory, there exists $\psi\in C^{\infty}(\Lambda^1T^*P)$ 
and $\omega\in C^{\infty}(\Lambda^2T^*P)$, with $d\omega=0$, such that $\theta=*d\psi+\omega$.  (Here, $\omega$ is
the sum of an exact form and a harmonic form.)  Therefore, $d\theta=d\!*\!d\psi$.  Letting $\Theta=*d\psi+d\psi
\in C^{\infty}(\Lambda^2_+T^*P)$, we have that $d\theta=d\Theta$.
\end{proof}

\begin{prop}\label{dimkerprop}
Use the notation of 
Definitions \ref{Hplusdfn} and \ref{l-l+dfn}.  
If $\mu\in[\lambda_-,-2)$,
$\dim\mathcal{K}(\mu)=\dim\mathcal{K}(-2)=b^2_+(N)$. If $\nu\in(-2,\lambda_+]$,
$$\dim\mathcal{K}(\nu)=b^2_+(N)-b^0(\Sigma)+b^1(\Sigma)+b^0(N)-b^1(N)+b^3(N).$$
\end{prop}

\begin{proof} Let $\mu\in[\lambda_-,-2)$ and $\nu\in(-2,\lambda_+]$.
By Theorem \ref{wSobthm}(a), $L^p_{k+1,\,\mu}\hookrightarrow
L^p_{k+1,\,-2}\hookrightarrow L^p_{k+1,\,\nu}$, so $\mathcal{K}(\mu)\subseteq\mathcal{K}(-2)\subseteq\mathcal{K}(\nu)$.

Recall the map $(d+d^*)_{\infty}$ defined by \eq{dd*inftyeq2} and the bundle $V_{\Sigma}$
given in \eq{VSigmaeq}.
By Definition \ref{jumpdfn}, $\d(-2)$ is the
dimension of the vector space of solutions $\xi(t,\sigma)$ to $(d+d^*)_{\infty}\xi=0$ of the form
$$\xi(t,\sigma)=e^{-2t}p(t,\sigma),$$
where $(t,\sigma)\in (T,\infty)\times\Sigma$ and $p(t,\sigma)$ is a polynomial in $t$ 
taking values in $C^{\infty}(V_{\Sigma}\otimes\C)$.  
Furthermore, by Proposition \ref{dmuprop} and its proof, $\d(-2)=b^1(\Sigma)$ and $\xi$ must be of the form
$$\xi(t,\sigma)=e^{-2t}\big(\alpha(\sigma)+dt\w *_{\Sigma}\alpha(\sigma)\big),$$
where $\alpha$ is a 
closed and coclosed 2-form on $\Sigma$ and $*_{\Sigma}$ is the Hodge star on $\Sigma$.   Suppose $\alpha\neq 0$. 
We shall see that the forms $\xi$ of this type correspond to forms on $N$ which add to the kernel of \eq{dd*eq} as
the rate crosses $-2$.

First, however, we must transform from the cylindrical picture to the conical one.  
Remember, in \eq{dd*inftyeq}, we related $(d+d^*)_{\infty}$ to the usual operators $d$ and $d^*$ acting on $m$-forms 
on $N\setminus K\cong (R,\infty)\times\Sigma$.
From this equation, since $\xi$ is a 2-form, 
we notice that we need to consider $\eta=e^{2t}\xi$ as the form on the cone related to $\xi$. 
Changing coordinates $(t,\sigma)$ to conical
coordinates $(r,\sigma)$, where $r=e^t$, we may write
\begin{equation}\label{etaeq}
\eta(r,\sigma)=\alpha(\sigma)+r^{-1}dr\w\ast_{\Sigma}\alpha(\sigma).
\end{equation}

We must now consider possible extensions of $\eta$ to a form on $N$.
Recall that $\Psi:(R,\infty)\times\Sigma\cong N\setminus K$ and that $\rho$ is a radius function on $N$ as in
Definition \ref{radiusfndfn}.  We also remind the reader of Proposition \ref{Upsilonprop}: that there is an isomorphism 
$$\tilde{\Upsilon}:\Lambda^2_+T^*\big((R,\infty)\times\Sigma\big)\rightarrow \Lambda^2_+T^*(N\setminus K).$$  
By the work in \cite[$\S$5]{LockhartMcOwen}, as $[\lambda_-,\lambda_+]\cap
\mathcal{D}=\{-2\}$, $\gamma\in\mathcal{K}(\nu)$ if and only if there exist $\eta$ of the form \eq{etaeq}, and
 $\zeta\in \mathcal{K}(\nu)$ with $\zeta$ asymptotic to $\tilde{\Upsilon}(\eta)$, such that
$\gamma-\zeta\in\mathcal{K}(\mu)$.  We can calculate:
$$|\eta(r,\sigma)|_{g_{\text{cone}}}=O(r^{-2})\qquad\text{as $r\rightarrow\infty$,}$$
where $g_{\text{cone}}$ is the conical metric on $(R,\infty)\times\Sigma$ as given in Definition \ref{ACsubmflddfn}, since $\alpha$
is a 2-form independent of $r$.  Thus, a form $\zeta$ asymptotic to $\tilde{\Upsilon}(\eta)$ must satisfy 
$\zeta=O(\rho^{-2})$ as $\rho\rightarrow\infty$, so 
$\zeta\notin\mathcal{K}(-2)$. Hence, $\gamma\notin\mathcal{K}(-2)$. 
 
Use the notation of Proposition \ref{seqprop}.
We show that 
a form $\zeta$ as in the previous paragraph will exist if and only if $[\alpha]\in p_2\big(H_{\dR}^2(N)\big)$.  
Necessity is clear, so we only consider sufficiency.  
Let $[\alpha]\in p_2\big(H_{\dR}^2(N)\big)$ with $\alpha$ 
the unique harmonic representative in the class.  Note that the map 
$\pi_{\Sigma}:(R,\infty)\times\Sigma\cong N\setminus K\rightarrow\Sigma$ induces an isomorphism
 $q_2:H^2_{\dR}(\Sigma)\rightarrow H^2_{\dR}(N\setminus K)$.  Now, $\eta$, given by
 \eq{etaeq}, is self-dual with respect to $g_{\text{cone}}$, so 
$$\zeta^{\prime}=\tilde{\Upsilon}(\eta)\in C^{\infty}\big(\Lambda^2_+T^*(N\setminus K)\big).$$
Moreover, since
$$\eta=\alpha+d(*_{\Sigma}\alpha\log r),$$
we see that $[\zeta^{\prime}]=[\eta]=[\alpha]$ in $H^2_{\dR}\big((R,\infty)\times\Sigma\big)$.  
Thus, $q_2^{-1}([\zeta^{\prime}])\in p_2\big(H^2_{\dR}(N)\big)$.  By the long exact sequence \eq{seqeq}, the image of
$p_2$ is equal to the kernel of $\partial_2$.  Hence, if $f\in C^{\infty}(N)$ such that $f=0$ on $K$ and $f=1$ on 
$\Psi\big((R+1,\infty)\times\Sigma\big)$, $[d(f\zeta^{\prime})]=0$ in $H^3_{\cs}(N)$.  Therefore, there exists a smooth,
 compactly supported, \emph{self-dual} 2-form $\chi$ on $N$ such that $d(f\zeta^{\prime})=d\chi$ by Lemma \ref{dimagelem}.
 
We define $\zeta=f\zeta^{\prime}-\chi$, which is a closed self-dual 2-form on $N$ asymptotic to $\tilde{\Upsilon}(\eta)$, since
it equals $\tilde{\Upsilon}(\eta)$ outside the support of $1-f$ and the compactly supported form $\chi$.  We have therefore
shown that the space of kernel forms added at rate $-2$ is isomorphic to the image of $p_2$.  We can calculate the dimension of this
image by \eq{seqeq}:
$$\text{dim $p_2\big(H_{\text{dR}}^2(N)\big)$}=-b^0(\Sigma)+b^1(\Sigma)+b^0(N)-b^1(N)+b^3(N).$$

By Proposition \ref{nochangeprop}, there are no changes in
$\mathcal{K}(\mu)$ for $\mu\in[\lambda_-,-2)$ and the argument
thus far shows that the kernel forms added at rate $-2$ do not
lie in $\mathcal{K}(-2)$.  We conclude that
$\mathcal{K}(\mu)=\mathcal{K}(-2)$ for $\mu\in[\lambda_-,-2)$, so $\dim\mathcal{K}(\mu)=
\dim\mathcal{K}(-2)=b^2_+(N)$ by Lemma \ref{dimkerlem}.
The latter part of the proposition follows from the formula and
arguments above, since $\mathcal{K}(\nu)$ does not alter for
$\nu\in(-2,\lambda_+]$ by Proposition \ref{nochangeprop}.
\end{proof}

\begin{note} By \eq{seqeq}, $-b^0(\Sigma)+b^1(\Sigma)+b^0(N)-b^1(N)+b^3(N)=0$
if $b^1(\Sigma)=0$.  
\end{note}

 Proposition \ref{dimkerprop} shows that
the function $k(\mu)=\dim\mathcal{K}(\mu)$ is
lower semi-continuous at $-2$ and is continuous there if 
$-2\notin\mathcal{D}$ by the note above. Our next result shows that 
$c(\mu)=\dim\mathcal{C}_+(\mu)$ is upper
semi-continuous at $-2$.

\begin{prop}\label{dimcokerprop2} 
Use the notation of Definition \ref{l-l+dfn}.  
If $\nu\in(-2,\lambda_+]$, then  $\dim\mathcal{C}_+(\nu)=\dim\mathcal{C}_+(-2)=b^3(N)$.  If
$\mu\in[\lambda_-,-2)$,
$$\dim\mathcal{C}_+(\mu)=b^0(\Sigma)-b^0(N)+b^1(N).$$
\end{prop}

\begin{proof}
We 
know the index of \eq{dd*eq} at 
rate $\nu$ by Propositions \ref{dimcokerprop} and \ref{dimkerprop}:
$$i_{\nu}(d+d^*)=b^2_+(N)-b^0(\Sigma)+b^1(\Sigma)+b^0(N)-b^1(N).$$
By Proposition \ref{jumpprop}, recalling that $\d(-2)=b^1(\Sigma)$, 
$$i_{\mu}(d+d^*)=b^2_+(N)-b^0(\Sigma)+b^0(N)-b^1(N).$$
Using Proposition \ref{conefunctionprop} again we get the dimension of $\mathcal{C}_+(\mu)$ claimed.

The only part left to do is to show that $\dim\mathcal{C}_+(-2)=b^3(N)$.  Remember that if the rate 
for the kernel forms is $\kappa$ then the `dual' rate for the cokernel forms is $-4-(\kappa-1)=-3-\kappa$.
Thus, by the work in \cite{LockhartMcOwen}, 
the forms $\gamma$ which are added to the cokernel as the rate crosses $-2$ from above must be asymptotic to homogeneous 
forms on the cone of order $O(r^{-1})$.  Therefore, $\gamma\notin\mathcal{C}(-2)\subseteq L^q_{l+1,\,-1}(\Lambda^3T^*N)$ as
required.
\end{proof}

\begin{note}
The difference $\dim\mathcal{C}_+(\mu)-\dim\mathcal{C}_+(\nu)$, as given in Proposition \ref{dimcokerprop2}, 
is equal to $\dim p_1\big(H^1_{\dR}(N)\big)$
as defined in Proposition \ref{seqprop}.   We can prove an isomorphism between the image of $p_1$ and the
cokernel forms added at $-2$ in a similar manner to the proof of Proposition \ref{dimkerprop}, but we omit it.  
\end{note}

Proposition \ref{jumpprop}, along with Propositions \ref{dimcokerprop}, 
\ref{dimkerprop} 
and \ref{dimcokerprop2}, give us the dimension of
$\mathcal{K}(\lambda)$, hence $\mathcal{M}(N,\lambda)$, for all
$\lambda\in [-2,0)\setminus\mathcal{D}$.  This is because the dimension of the cokernel is 
constant for these choices of $\lambda$, so the change of index of \eq{dd*eq} equals the change in the
dimension of $\mathcal{K}(\lambda)$.

\begin{prop}\label{dimprop1}
In the notation of Definitions 
\ref{jumpdfn} and \ref{Hplusdfn},
 if $\lambda\in [-2,0)\setminus\mathcal{D}$, 
$$\dim\mathcal{M}(N,\lambda)=b_+^2(N) -b^0(\Sigma)+b^1(\Sigma)+b^0(N)-b^1(N)+b^3(N)
+\!\!\!\!\!\!\sum_{\mu\in\mathcal{D}\,\cap(-2,\,\lambda)}\!\!\!\!\!\!\d(\mu).$$
\end{prop}

\subsection{Rates {\boldmath $\lambda\in [0,1)\setminus\mathcal{D}$}}

We now discuss the case $\lambda\geq0$ and begin by studying the point
$0\in\mathcal{D}$.  
Recall that $\d(0)$, as given by
Proposition \ref{jumpdfn}, is equal to the dimension of
$$D(0)=\{(\alpha,\beta)\in C^{\infty}(\Lambda^2T^*\Sigma\oplus
\Lambda^3T^*\Sigma)\,:\,d_{\Sigma}\alpha=0\;\text{and}\;d_{\Sigma}\!*_{\Sigma}\!\alpha+d^*_{\Sigma}\beta=2\alpha\},$$
using the notation of Proposition \ref{Dprop}.

It is clear, using integration by parts, that the equations which
$(\alpha,\beta)\in D(0)$ satisfy are equivalent to
\begin{equation*}\label{d0eq}
d_{\Sigma}\!*_{\Sigma}\!\alpha=2\alpha\quad\text{and}\quad d^*_{\Sigma}\beta=0.\end{equation*}
The latter equation corresponds to constant 3-forms on $\Sigma$, so
the solution set has dimension equal to $b^0(\Sigma)$. If we define
\begin{equation}\label{Zeq}
Z=\{\alpha\in
C^{\infty}(\Lambda^2T^*\Sigma)\,:\,d_{\Sigma}\!*_{\Sigma}\!\alpha=2\alpha\},
\end{equation}
then $\d(0)=b^0(\Sigma)+\text{dim}\, Z$.

\medskip

Suppose that $\beta\in C^{\infty}(\Lambda^3T^*\Sigma)$ satisfies
$d^*_{\Sigma}\beta=0$ and corresponds to a form on $N$ which adds to the
kernel of \eq{dd*eq} at rate $0$.  There exists, by \cite[$\S$5]{LockhartMcOwen}, a 4-form $\zeta$ on $N$ asymptotic to
 $r^3dr\w\beta$ on $N\setminus K\cong
(R,\infty)\times\Sigma$, and a self-dual 2-form $\xi$ of order
$o(1)$ as $\rho\rightarrow\infty$ such that
$d\xi+d^*\zeta=0$, where $\rho$ is a radius function on $N$.

Since $d^*\zeta$ is exact, $*\zeta$ is a harmonic function which
is asymptotic to a function $c$, constant on each \emph{end} of $N$,
as given in Definition \ref{ACmflddfn}. Applying \cite[Theorem
7.10]{Joyce5} gives a unique harmonic function $f$ on $N$ which
converges to $c$ with order $O(\rho^{\mu})$ for all $\mu\in(-2,0)$.
The theorem cited is stated for an AC special Lagrangian submanifold
$L$, but only uses the fact that $L$ is an AC manifold
and hence is applicable here. Therefore, $*\zeta-f=o(1)$ as
$\rho\rightarrow\infty$ and hence, by the Maximum Principle,
$*\zeta=f$.

We deduce that $d^*\zeta$ and $d\xi$ are
$O(\rho^{-3+\epsilon})$ as $\rho\rightarrow\infty$ for any
$\epsilon>0$ small, hence they lie in $L^2$.  Integration by parts,
now justified, shows that $d^*\zeta=0$ and we conclude that
$*\zeta$ is constant on each component of $N$.
 Hence, the piece of $b^0(\Sigma)$ in $\d(0)$ that adds
to the dimension of the kernel is equal to $b^0(N)$.  The
other 3-forms $\beta$ on $\Sigma$ satisfying $d^*_{\Sigma}\beta=0$ must
correspond to cokernel forms, so $b^0(\Sigma)-b^0(N)$ is
subtracted from the dimension of the cokernel at $0$.

\medskip

For each end of $N$, we can define self-dual 2-forms of order $O(1)$
given by translations of it, written as $\frac{\partial}{\partial
x_j}\cdot\varphi$ for $j=1,\ldots,7$.  If the end is a flat
$\R^4$, we only get three such self-dual 2-forms from it. So, if
$k^{\prime}$ is the number of ends which are not 4-planes,
$$\text{dim}\,Z\geq
7k^{\prime}+3(b^0(\Sigma)-k^{\prime})=3b^0(\Sigma)+4k^{\prime}.$$
Moreover, the translations of the components of $N$ must correspond
to kernel forms, since they are genuine deformations of $N$.  If a component of $N$ is a 4-plane then there are only 
three nontrivial translations of it.  
Therefore, if $k$ is the number of components of $N$ which are not
4-planes, at least $3b^0(N)+4k$ is added to the dimension of the
kernel at rate $0$ from $\text{dim}\, Z$.

\medskip

From this discussion, we can now state and prove the following inequalities for the dimension of
$\mathcal{M}(N,\lambda)$ for $\lambda\in[0,1)\setminus\mathcal{D}$.

\begin{prop}\label{dimprop2}
Use the notation from Proposition \ref{jumpdfn}, Definition \ref{Hplusdfn} and Theorem \ref{modulispacethm}. 
If $\lambda\in [0,1)\setminus\mathcal{D}$,
$$\dim\mathcal{M}(N,\lambda)\leq\dim\mathcal{K}(0)+b^0(N) +\dim 
Z-\dim\mathcal{B}(\lambda)+\!\!\!\!\!\!\sum_{\mu\in\mathcal{D}\,
\cap(0,\,\lambda)}\!\!\!\!\!\!\d(\mu)$$ where $Z$ is defined in \eq{Zeq} and
$$\dim\mathcal{K}(0)=b_+^2(N)
-b^0(\Sigma)+b^1(\Sigma)+b^0(N)-b^1(N)+b^3(N)
+\!\!\!\!\!\!\sum_{\mu\in\mathcal{D}\,\cap(-2,\,0)}\!\!\!\!\!\!\d(\mu).$$
Moreover, if $k$ is the number of components of $N$ which are not a
flat $\R^4$,
$$\dim\mathcal{M}(N,\lambda)\geq
\dim\mathcal{K}(0)+b^0(\Sigma)+3b^0(N)+4k-b^3(N)
-\dim\mathcal{B}(\lambda)+
\!\!\!\!\!\!\sum_{\mu\in\mathcal{D}\,
\cap(0,\,\lambda)}\!\!\!\!\!\!\d(\mu).$$
\end{prop}

\begin{proof}
By Theorem \ref{modulispacethm}, $\dim\mathcal{M}(N,\lambda)=\dim\mathcal{K}(\lambda)-\dim\mathcal{B}(\lambda)$.
Using Propositions \ref{jumpprop}, \ref{dimcokerprop} and \ref{dimkerprop} and the notation of Definition \ref{l-l+dfn}, we calculate: 
\begin{align*}
i_{\lambda}(d+d^*)&=i_{\lambda_+}(d+d^*)+\!\!\!\!\!\!\sum_{\mu\in\mathcal{D}\,\cap(-2,\,\lambda)}\!\!\!\!\!\!\d(\mu) \\
&=b^2_+(N)-b^0(\Sigma)+b^1(\Sigma)+b^0(N)-b^1(N)+\!\!\!\!\!\!\sum_{\mu\in\mathcal{D}\,\cap(-2,\,\lambda)}\!\!\!\!\!\!\d(\mu)\\
&=\dim\mathcal{K}(0)-b^3(N)+\d(0)+\!\!\!\!\!\!\sum_{\mu\in\mathcal{D}\,\cap(0,\,\lambda)}\!\!\!\!\!\!\d(\mu).
\end{align*}
(An argument using asymptotic behaviour of kernel forms as in the proof of Proposition \ref{dimkerprop} leads to 
lower semi-continuity of $k(\mu)=\dim\mathcal{K}(\mu)$ at $0$).  
From the discussion above, the part of $\d(0)=b^0(\Sigma)+\text{dim}\,Z$ 
that subtracts from the cokernel of \eq{dd*eq} as the rate crosses $0$ is greater than or equal to
 $b^0(\Sigma)-b^0(N)$.
Therefore, $\dim\mathcal{C}_+(\lambda)\leq b^3(N)-b^0(\Sigma)+b^0(N)$ using Proposition \ref{dimcokerprop}, so 
\begin{align*}
\dim\mathcal{K}(\lambda)&=i_{\lambda}(d+d^*)+\dim\mathcal{C}_+(\lambda)\\
&\leq  i_{\lambda}(d+d^*)+b^3(N)-b^0(\Sigma)+b^0(N)\\
&=\dim\mathcal{K}(0)+b^0(N)+\dim Z+\!\!\!\!\!\!\sum_{\mu\in\mathcal{D}\,\cap(0,\,\lambda)}\!\!\!\!\!\!\d(\mu).
\end{align*}
The first inequality follows.

Now, from the discussion above, 
at least $b^0(N)+3b^0(N)+4k$ is added to the kernel of \eq{dd*eq}
 from $\d(0)$ as
the rate crosses $0$.  
Remember $\dim\mathcal{C}_+(\lambda)\leq b^3(N)-b^0(\Sigma)+b^0(N)$.  Therefore, at most this amount could be added to the
index of \eq{dd*eq} from the sum of $\d(\mu)$ terms \emph{without} adding to the kernel.  Thus,
\begin{align*}
\dim\mathcal{K}(\lambda)&\geq\dim\mathcal{K}(0)+4b^0(N)+4k-b^3(N)+\!\!\!\!\!\!\sum_{\mu\in\mathcal{D}\,\cap(0,\,\lambda)}\!\!\!\!\!\!\d(\mu)\\
&\quad-b^3(N)+b^0(\Sigma)-b^0(N)\\
&=\dim\mathcal{K}(0)+b^0(\Sigma)+3b^0(N)+4k-b^3(N)+\!\!\!\!\!\!\sum_{\mu\in\mathcal{D}\,\cap(0,\,\lambda)}\!\!\!\!\!\!\d(\mu),
\end{align*}
from which we can easily deduce the second inequality.
\end{proof}

\subsection{Rates {\boldmath $\lambda<-2$, $\lambda\notin\mathcal{D}$}}

By Theorem \ref{modulispacethm} and its proof, we can think of the map $\pi$ as a projection from an open neighbourhood of $0$ 
in the infinitesimal deformation space $\mathcal{I}(N,\lambda)$ to the obstruction space $\mathcal{O}(N,\lambda)$.  Thus, 
the expected dimension of the moduli space is $\dim\mathcal{I}(N,\lambda)-\dim
\mathcal{O}(N,\lambda)$.  We now find a lower bound for this dimension.

\begin{prop}\label{dimprop3}
Use the notation of 
Definitions \ref{infdefspacedfn}, \ref{obsspacedfn} and \ref{Hplusdfn}.  The expected dimension of the moduli space satisfies
\begin{align*}
\dim\mathcal{I}(N,\lambda)&-\dim\mathcal{O}(N,\lambda)\geq \\
&
b^2_+(N)-b^0(\Sigma)+b^0(N)-b^1(N)+b^3(N)-\!\!\!\!\!\!\sum_{\mu\in\mathcal{D}\,
\cap(\lambda,\,-2)}\!\!\!\!\!\!\d(\mu).
\end{align*}
\end{prop}

\begin{proof}
In Definition \ref{obsspacedfn} we noted that
$$\dim\mathcal{O}(N,\lambda)=\dim\mathcal{C}_+(\lambda)-\dim\mathcal{C}(\lambda),$$
where $\mathcal{C}(\lambda)$ is given in \eq{fullcokereq}.  Since $L^q_{l+1,\,-1}\hookrightarrow L^q_{l+1,\,-3-\lambda}$ 
 by 
Theorem \ref{wSobthm}(a) (as $\lambda>-2$), we deduce that $\dim\mathcal{C}(\lambda)\geq b^3(N)$ by Proposition \ref{dimcokerprop}.  Therefore,  
$$\dim\mathcal{I}(N,\lambda)-\dim\mathcal{O}(N,\lambda)\geq \dim\mathcal{K}(\lambda)-\dim\mathcal{C}_+(\lambda)+b^3(N).$$
Applying Propositions \ref{jumpprop}, \ref{dimkerprop} and \ref{dimcokerprop2} completes the proof.
\end{proof}

\begin{remark}
Any nonplanar coassociative 4-fold 
which is AC with rate $\lambda$ has a 1-dimensional family of 
nontrivial AC coassociative deformations given by dilations.  The obvious question is: where does this appear in the calculation
of the dimension of the moduli space?  The author believes that it should come from kernel forms 
of rate $\mu$, where $\mu$ is the greatest element of $\mathcal{D}\cap(-\infty,\lambda)$.
\end{remark}

\section{Invariants}
\label{invsection} 

Suppose for this section that $N$ is a coassociative 4-fold which is
AC to a cone $C$ in $\R^7$ with rate $\lambda<-2$ and suppose, for convenience, that $N$ is
connected.  Recall the space $\mathcal{H}^2_+(N)$ and topological invariant $b^2_+(N)$ given in Definition \ref{Hplusdfn}.

\medskip

Consider the deformation $N\mapsto e^t N$ for $t\in\R$.  Clearly $e^tN$ is
coassociative and AC to $C$ with rate $\lambda$ for all $t\in\R$. Let $u$ be
the dilation vector field on $\R^7$ as given in \eq{udilationeq}.
 Define a self-dual 2-form $\alpha_u$ on $N$, as in the note after Proposition \ref{jmathprop}, by
\begin{equation}\label{alphaueq}
\alpha_u=(u\cdot\varphi)|_{TN}.
\end{equation}  
The deformation
corresponding to $t\alpha_u$, as given in Definition \ref{Fdfn},  is $e^tN$.  

We calculate, using \eq{uLiederiveq},
$$d\alpha_u=d(u\cdot\varphi)|_{TN} =3\varphi|_{TN} =0$$
since $N$ is coassociative.  As $\lambda<-2$, $\alpha_u\in
L^2(\Lambda^2_+T^*N)$ and hence $t\alpha_u\in\mathcal{H}_+^2(N)$ for all $t\in\R$.  So, there is a 1-dimensional 
subspace of $\mathcal{H}_+^2(N)$ whenever $e^tN\neq N$ for some $t\neq 0$; that is, when $N$ is not
a cone.
Consequently, $b_+^2(N)=0$ forces $N\cong\R^4$, as $N$ is nonsingular.

This discussion leads to our first invariant.

\begin{dfn}\label{XNdfn} Let $N$ be an AC coassociative 4-fold in $\R^7$ with rate $\lambda<-2$, define $\alpha_u$ by \eq{alphaueq} and let $X(N)=\|\alpha_u\|_{L^2}^2$.  
By the discussion above, this 
 is a well-defined invariant of $N$.
\end{dfn}

We now consider, for $\Gamma$ a 2-cycle in $N$ and $D$ a 3-cycle in $\R^7$ such that $\partial D =\Gamma$, the integral
$\int_D\varphi$.  

Suppose $D^\prime$ is another 3-cycle with $\partial D^\prime=\Gamma$. Then
$\partial(D-D^\prime)=0$, so
$$\int_{D-D^\prime}\varphi=[\varphi]\cdot[D-D^\prime]=0$$
since $[\varphi]\in H_{\dR}^3(\R^7)$ is zero.  Therefore
$\int_D\varphi=\int_{D^\prime}\varphi$; that is, the integral is independent of the choice of 3-cycle with boundary $\Gamma$.

Suppose instead that
$\Gamma^\prime$ is a 2-cycle in $N$ such that
$\Gamma-\Gamma^\prime=\partial E$ for some 3-cycle $E\subseteq N$.
Let $D$ and $D^\prime$ be 3-cycles such that $\partial D=\Gamma$ and
$\partial D^\prime=\Gamma^\prime$.  Then $\partial D=\partial
(E+D^\prime)=\Gamma$ and thus
$$\int_D\varphi
=\int_{E+D^\prime}\varphi
=\int_E\varphi+\int_{D^\prime}\varphi =\int_{D^\prime}\varphi,
$$
since $E\subseteq N$ and $\varphi|_N=0$.  Hence, we have shown the integral only depends on the homology class $[\Gamma]$.

We can now define our second invariant for any coassociative 4-fold in $\R^7$.

\begin{dfn}\label{YNdfn}  Let $N$ be a coassociative 4-fold in $\R^7$.
Let $\Gamma$ be a 2-cycle in $N$ and let $D$
be a 3-cycle in $\R^7$ such that $\partial D=\Gamma$.  Define
$[Y(N)]\in H^2_{\dR}(N)$ by:
$$[Y(N)]\cdot[\Gamma]=\int_D\varphi.$$
The work above shows that this is well-defined. 
\end{dfn}

\noindent 
In the notation of Definition \ref{YNdfn}, 
we calculate, using \eq{uLiederiveq} and \eq{alphaueq},
$$3\int_D\varphi=\int_Dd(u\cdot\varphi)=\int_{\Gamma}u\cdot\varphi=\int_{\Gamma}\alpha_u=[\alpha_u]\cdot[\Gamma].$$
Hence, $3[Y(N)]=[\alpha_u]$.

We remind the reader of Definitions \ref{Hmdfn} and \ref{Hplusdfn}, 
in particular that $\mathcal{J}(N)$ is the image of $H^2_{\cs}(N)$ in
$H^2_{\dR}(N)$.  We showed that $\alpha_u\in\mathcal{H}^2_+(N)$ and we know that 
$\mathcal{H}^2(N)\cong\mathcal{J}(N)$ via $\gamma\mapsto[\gamma]$. 
 Thus, $[Y(N)]$ and $[\alpha_u]$ lie in $\mathcal{J}(N)$.  
Using the product \eq{producteq}, we deduce that
$$9[Y(N)]^2=[\alpha_u]\cup[\alpha_u]=\int_N\alpha_u\w\alpha_u =
\|\alpha_u\|_{L^2}^2 =X(N).$$
We have used the fact that if $\alpha_u=\eta+d\xi$, for some compactly
supported closed 2-form $\eta$ and $\xi\in
C^{\infty}_{\lambda+1}(\Lambda^1T^*N)$, then an integration by parts
argument, valid since $\lambda<-2$, shows that
$$\int_N\alpha_u\w\alpha_u=\int_N\eta\w\eta.$$

\medskip

We may thus derive a test which determines whether $N$ is a cone.

\begin{prop}\label{Nconeprop}
 Let $N$ be a connected coassociative
4-fold in $\R^7$ which is AC with rate $\lambda<-2$.  Use the notation of Definitions \ref{Hplusdfn}, 
\ref{XNdfn} and \ref{YNdfn}.  If
$b_+^2(N)=0$, $N$ is a cone and hence
a linear $\R^4$ in $\R^7$ as $N$ is nonsingular.  Moreover, $N$ is a cone if and only if
$X(N)=0$ or, equivalently, $[Y(N)]^2=0$.
\end{prop}

We note, for interest, that a similar argument holds for a connected
special Lagrangian $m$-fold $L$ in $\C^m$ ($m\geq 2$) which is AC to a cone $C\cong(0,\infty)\times\Sigma$ with rate $\lambda<-m/2$, in the
following sense. 
Let $u$ be the dilation vector field on $\C^m$, and let $\omega_m$ and $\Omega_m$ be the symplectic and holomorphic volume forms
on $\C^m$ respectively.  Define 
$$\beta_u=(u\cdot\omega_m)|_{TL}\quad\text{and}\quad \gamma_u=(u\cdot\Omega_m)|_{TL},$$
which are both zero if and only if $L$ is a cone.  
By the proof of \cite[Theorem 3.6]{McLean}, $\beta_u$ and $\gamma_u$ are closed and $\gamma_u=*\beta_u$.  
Moreover, by \eq{Lpeq} and Theorem \ref{wSobthm}(a), $\beta_u$ and $*\beta_u$ lie in $L^2$ as $\lambda<-m/2$.

By \cite[Example 0.15]{Lockhart}, the closed and
coclosed
 $(m-1)$-forms in $L^2$ uniquely represent the cohomology classes in 
$H^{m-1}_{\dR}(L)$.
Therefore, if $b^{m-1}(L)=0$, $L$ is a cone and hence a linear
$\R^m$ in $\C^m$ as $L$ is nonsingular.
 Moreover, we may define the invariant $X(L)=\|\beta_u\|_{L^2}^2$ in an
analogous way to $X(N)$ and see that $L$ is a cone if and only if
$X(L)=0$.

We define two invariants, $[Y(L)]\!\in\! H^1_{\dR}(L)$ and
$[Z(L)]\!\in\! H^{m-1}_{\dR}(L)$ by
\begin{align*}
[Y(L)]\cdot[\Gamma]=\int_D\omega_m\quad\text{and}
\quad[Z(L)]\cdot[\Gamma^{\prime}]=\int_{D^{\prime}}\Im\Omega_m,
\end{align*}  where
$\Gamma$ is a 1-cycle in $L$, $D$ is a 2-cycle in $\C^m$ such that $\partial
D=\Gamma$, $\Gamma^{\prime}$ is an $(m-1)$-cycle in $L$ and
$D^{\prime}$ is an $m$-cycle in $\C^m$ such that $\partial
D^{\prime}=\Gamma^{\prime}$.  Analogues of these invariants on $\Sigma$ 
were introduced in \cite[Definition 7.2]{Joyce5}.  

Using similar calculations from our discussion of $[Y(N)]$, we can show
that $2[Y(L)]=[\beta_u]$ and $m[Z(L)]=[*\beta_u]$. 
By \cite[Proposition 7.3]{Joyce5}, $[Y(L)]$ lies 
in the image of the compactly supported cohomology for rates less than $-1$.
Therefore, for $\lambda<-m/2\leq-1$ (as $m\geq 2$), we calculate:
$$2m[Y(L)]\cup[Z(L)]=[\beta_u]\cup[*\beta_u]
=\int_L\beta_u\w*\beta_u=\|\beta_u\|_{L^2}^2=X(L),$$
where the cup product is given by the pairing between $H^1_{\cs}(L)$ and $H^m_{\dR}(L)$, since $[\beta_u]$
can be represented by a compactly supported closed 1-form.
Hence $X(L)=0$ if and only if $[Y(L)]\cup[Z(L)]=0$. 

We write these observations as a
proposition.

\begin{prop}\label{SLconeprop} Let $L$ be a connected special Lagrangian $m$-fold
in $\C^m$ ($m\geq 2$) which is AC with rate $\lambda<-m/2$.  If $b^{m-1}(L)=0$, 
$L$ is a cone and hence a linear $\R^m$ in $\C^m$ as $L$ is nonsingular.  Moreover, in the notation above,
$L$ is a cone if and only if $X(L)=0$ or, equivalently,
$[Y(L)]\cup [Z(L)]=0$.
\end{prop}

\section{Examples}
\label{exsection}

Before discussing examples we prove a general result which shows that if an AC coassociative
4-fold converges sufficiently quickly to a cone with symmetries, then it inherits those symmetries.  The
analogous statement is already known for AC special Lagrangian submanifolds.

\begin{prop}\label{conesymprop}
Let $N$ be a coassociative 4-fold which is asymptotically conical to a cone $C$ in $\R^7$ with rate $\lambda<-2$.  
Let $\GG$ be the identity component of the automorphisms of $C$, $\Aut(C)\subseteq 
\GG_2$.  Then $N$ is $\GG$-invariant.
\end{prop}

\begin{proof}
Suppose for a contradiction that $N$ is not $\GG$-invariant.  Then there exists a vector field $v$ in the 
Lie algebra $\mathfrak{g}$ of $\GG$ which preserves $C$ but not $N$.  
Since $v$ does not preserve $N$, it defines a nonzero closed self-dual 2-form 
$\alpha_v=(v\cdot\varphi)|_{TN}$ as in the note after Proposition \ref{jmathprop}.  
Morever, $\alpha_v\in C^{\infty}_{\lambda}(\Lambda^2_+T^*N)$
since the deformation of $N$ defined by $\alpha_v$, as in Definition \ref{Fdfn}, converges to $C$ at the same rate as $N$,
recalling that $v$ preserves $C$.  

 Notice that, since $v\in\mathfrak{g}$, $\GG\subseteq\GG_2$ and $\varphi$ is invariant under $\GG_2$, the
class of $\varphi$ in $H^3(\R^7,N)$, which is zero, is unchanged under flow along $v$.  
Let $u$ be the dilation vector field on $\R^7$ as given in \eq{udilationeq}. 
This leads us to define $\beta_v\in C^{\infty}(\Lambda^1T^*N)$ by
$$\beta_v(x)=\frac{1}{3}\,\varphi(v,u,x)$$
for $x\in TN$.  Thus, by \eq{uLiederiveq}, $d\beta_v=\alpha_v$ and, since $\alpha_v\in C^{\infty}_{\lambda}(\Lambda^2_+T^*N)$ and $u=O(r)$ as $r\rightarrow\infty$, 
where $r$ is the radial distance 
in $\R^7$, $\beta_v\in C^{\infty}_{\lambda+1}(\Lambda^1T^*N)$.  

Since $\alpha_v=d\beta_v$ is closed and self-dual, $d^*d\beta_v=0$.  If $\rho$ is a radius function on $N$ as in Definition
\ref{radiusfndfn}, $d\beta_v=O(\rho^{\lambda})$ and $\beta_v=O(\rho^{\lambda+1})$ as $\rho\rightarrow\infty$.  As 
$L^2_{0,\,-2}=L^2$ by \eq{Lpeq} and $\lambda<-2$, we see from Theorem \ref{wSobthm}(a) that $d\beta_v\in L^2$ and the following integration by parts argument
is valid:
$$\|\alpha_v\|_{L^2}^2=\|d\beta_v\|_{L^2}^2=\langle d^*d\beta,\beta \rangle_{L^2}=0.$$
Therefore, $\alpha_v=0$, our required contradiction.
\end{proof}

\subsection{{\boldmath $\SU(2)$}-invariant AC coassociative 4-folds}
The coassociative 4-folds we discuss were constructed in
\cite[$\!$Theorem IV$\!$.3.2]{HarLaw}. The construction involves the
consideration of $\R^7$ as the imaginary part of the octonions,
$\O$, and the definition of an action of $\SU(2)$ on $\Im\O$ as
unit elements in the
 quaternions, $\H$.   We first define the action then state the result.

\begin{dfn}\label{su2actiondfn}
Write $\Im\O\cong\Im\H\oplus\H$.
Define an action on $\Im\O$ by
$$(x,y)\mapsto (qx\bar{q},\bar{q}y)$$
for $q\in\H$ such that $|q|=1$.  This defines an action of $\SU(2)$ on $\Im\O\cong\R^7$.
\end{dfn} 

\begin{prop}\label{su2exprop} 
Use the notation of Definition \ref{su2actiondfn}.  Let $(e,f)\in\Im\H\oplus\H\cong\Im\O$ be such that $|e|=|f|=1$ and let $c\in\R$.
Then
\begin{align*}M_c=\{(sqe\bar{q},r\bar{q}f)\,:\,q\in\H,\,|q|=1,\,
\text{and}\;\,
s(4s^2-5r^2)^2=c\;\text{for}\;r\geq0,\,s\in\R\}\end{align*} is an $\SU(2)$-invariant
coassociative 4-fold in $\R^7\cong\Im\O$, where the action is given in Definition \ref{su2actiondfn}.  Moreover,
every coassociative 4-fold in $\R^7$ invariant under this $\SU(2)$ action is of the form $M_c$ for some choice of $c$, $e$ and $f$.
\end{prop}

Suppose first that $c=0$.  So, $s=0$, $s=\frac{\sqrt{5}}{2}r$ or $s=-\frac{\sqrt{5}}{2}r$.  Let
\begin{align*}
M_0^0&=\{(0,r\bar{q}f)\,:\,r\geq0\}\cong\R^4
\qquad\text{and}\\[4pt] 
M_0^{\pm}&=\left\{r\left(\pm\frac{\sqrt{5}}{2}\,qe\bar{q},\bar{q}f\right)\,:\,r>0\right\}\cong(0,\infty)\times\mathcal{S}^3,
\end{align*}
so that $M_0=M_0^0\sqcup M_0^+\sqcup M_0^-$.  Thus, $M_0$ is a $\SU(2)$-invariant cone with three ends, two of which are
 diffeomorphic to cones on $\mathcal{S}^3$ and one which is a flat $\R^4$.

\medskip

Suppose now $c\neq 0$ and take $c>0$ without loss of generality.  This
forces $s>0$ and $4s^2-5r^2\neq 0$.  It is then clear that $M_c$ has
two components, $M_c^+$ and $M_c^-$, corresponding to $4s^2-5r^2>0$
and $4s^2-5r^2<0$ respectively.

\smallskip

The component $M_c^+$ is AC 
to the cone $M_0^+$.  We calculate the rate as follows.
For large $r$, $s$ is approximately equal to $\frac{\sqrt{5}}{2}\,r$
and hence $4s^2-5r^2=O(r^{-\frac{1}{2}})$. Therefore
$s=\frac{\sqrt{5}}{2}\,r+O(r^{-\frac{3}{2}})$ and thus $M_c^+$
converges with rate $-3/2$ to $M_0^+$.  For
each $r\neq 0$ we have an $\mathcal{S}^3$ orbit in $M_c^+$, but when
$r=0$ there is an $\mathcal{S}^2$ orbit. Therefore, topologically,
$M_c^+$ is an $\R^2$ bundle over $\mathcal{S}^2$. Hence
$H^2_{\dR}(M_c^+)=\R$. 

Suppose, for a contradiction, that $b^2_+(M_c^+)=1$, in the notation of Definition
\ref{Hplusdfn}.
Therefore, there exists a smooth, closed, self-dual 2-form in $L^2=L^2_{0,\,-2}$.  The existence of this form, and 
the fact that the deformation theory of $M_c^+$ is unobstructed by Corollary \ref{modulispacecor}, 
means that there is a coassociative deformation $\tilde{M}_c^+$ of $M_c^+$  which is AC to $M_0^+$ with rate less than $-2$.  
However, $M_0^+$ is $\SU(2)$-invariant so, by 
Proposition \ref{conesymprop}, $\tilde{M}_c^+$ must itself be invariant under $\SU(2)$.  Proposition \ref{su2exprop} describes the
family of all coassociative 4-folds invariant under this action, so $\tilde{M}_c^+$ must be AC with rate $-3/2$, our
required contradiction. 

Hence, $b^2_+(M_c^+)=0$ and,
since the $\SU(2)$ action has generic orbit $\mathcal{S}^3$, we
conclude that $M_c^+$ is isomorphic to the bundle $\mathcal{O}(-1)$
over $\C\P^1$.

\smallskip

We now turn to $M_c^-$.  Here there are two ends, one of which has
the same behaviour as the end of $M_c^+$ and the other is where
$s\rightarrow 0$. For the latter case, we quickly see that
$s=O(r^{-4})$ and so $M_c^-$ converges at rate $-4$ to $M_0^0$. As
the case $r=0$ is excluded here, $M_c^-$ is topologically
$\R\times\mathcal{S}^3$ and converges with rate $-3/2$ to
$M_0^+$ and rate $-4$ to $M_0^0$ at its two
ends. Hence $H^2_{\dR}(M_c^-)=0$.

\smallskip

Consequently, for $c\neq 0$,
 $b_+^2(M_c)=0$.  Notice also that the link 
$\Sigma$ of the cone to which $M_c$ converges has $b^1(\Sigma)=0$.  
Applying Corollary \ref{modulispacecor} and Proposition \ref{dimprop1},  
the moduli space $\mathcal{M}(M_c,-3/2)$ is a smooth
manifold with dimension
$$\sum_{\mu\in\mathcal{D}\cap(-2,-3/2)}\!\!\!\!\!\!\d(\mu),$$
using the notation of Proposition \ref{Dprop} and Definition \ref{jumpdfn}.

\begin{note}
Fox \cite[Example 9.2 \& Theorem 9.3]{Fox} showed that the examples given in Proposition \ref{su2exprop} are special cases 
of a surface bundle construction over a \emph{pseudoholomorphic curve} in $\mathcal{S}^6$ with \emph{null-torsion} (in 
Proposition \ref{su2exprop} the 
curve is a totally geodesic 2-sphere).  The examples he produces are also coassociative 4-folds $N$ which are AC to some cone $C$ 
with rate $-3/2$, but the link $\Sigma$ of $C$ may now have $b^1(\Sigma)>0$ and $b^2_+(N)$ may be non-zero.
\end{note}

\subsection{2-ruled coassociative 4-folds}
In \cite{Lotay2}, the author introduced the notion of \emph{2-ruled}
4-folds.  A 4-dimensional submanifold $M$ in $\R^n$ is
\emph{2-ruled} if it admits a fibration over a 2-fold $S$ such that
each fibre is an affine 2-plane. As commented in
\cite[$\S$3]{Lotay2}, if $S$ is compact then $M$ is AC with rate $\lambda\leq0$ to its \emph{asymptotic cone}.  Therefore,
2-ruled 4-folds can be examples of AC submanifolds which are not
cones. Unfortunately, the examples in \cite[$\S$5]{Lotay2} of 2-ruled
4-folds are not coassociative, but are calibrated submanifolds of
$\R^8$ known as \emph{Cayley} 4-folds.  However, the author has derived a system of ordinary differential equations whose
solutions define 2-ruled coassociative 4-folds invariant under a $\U(1)$ action.  Though the solutions have not yet been
found, it is clear that this system will provide AC examples.

The result \cite[Theorem 3.5.1]{Bundle} gives a method of
construction for coassociative 4-folds in $\R^7$ which are
necessarily 2-ruled. However, the examples in \cite{Bundle} of
coassociative 4-folds are not AC in the notion of Definition
\ref{ACsubmflddfn} but rather in a much weaker sense.
In \cite{Fox}, Fox studies 2-ruled coassociative 4-folds, with a particular focus on
cones.  His works leads to a number of techniques for constructing
2-ruled coassociative 4-folds.  

Recently the author (\cite[Theorem 7.5]{Lotaylag}) has explicitly classified the 2-ruled coassociative cones and showed that the
 general cone can be constructed using a pseudoholomorphic curve in $\mathcal{S}^6$ and holomorphic
data on the curve.  It is the hope of the author to use this result to produce nontrivial
examples of AC 2-ruled coassociative 4-folds, using both the methods in \cite{Lotay2} and other means.

\noindent\textsl{Jason D.~Lotay\\ University College\\ High Street\\ Oxford\\ OX1 4BH\\ lotayj@maths.ox.ac.uk}


\begin{thebibliography}{99}
\bibitem{Bartnik} R.~Bartnik, {\it The Mass of an Asymptotically
Flat Manifold}, Comm.~Pure Appl.~Math. {\bf 39}, 661-693, 1986.
\bibitem{Fox} D.~Fox, {\it Coassociative Cones that are Ruled by 2-Planes}, Asian J.~Math. {\bf 11}, 535-554, 2007.
\bibitem{HarLaw} R.~Harvey and H.~B.~Lawson, {\it Calibrated Geometries},
Acta Math. {\bf 148}, 47-152, 1982.
\bibitem{Bundle} M.~Ionel, S.~Karigiannis and M.~Min-Oo,
{\it Bundle Constructions of Calibrated Submanifolds of $\R^7$ and
$\R^8$}, Math.~Res.~Lett. {\bf 12}, 493-512, 2005.
\bibitem{Joyce5} D.~D.~Joyce, {\it Special Lagrangian Submanifolds
with Isolated Conical Singularities. I.~Regularity}, Ann.~Global
Anal.~Geom. {\bf 25}, 201-251, 2004.
\bibitem{Joyce6} D.~D.~Joyce, {\it Special Lagrangian Submanifolds
with Isolated Conical Singularities. II.~Moduli Spaces}, Ann.~Global
Anal.~Geom. {\bf 25}, 301-352, 2004.
\bibitem{Joyce7} D.~D.~Joyce, {\it Special Lagrangian Submanifolds
with Isolated Conical Singularities. III.~Desingularization, The
Unobstructed Case}, Ann.~Global Anal.~Geom. {\bf 26}, 1-58, 2004.
\bibitem{Joyce8} D.~D.~Joyce, {\it Special Lagrangian Submanifolds
with Isolated Conical Singularities. IV.~Desingularization,
Obstructions and Families}, Ann.~Global Anal.~Geom. {\bf 26},
117-174, 2004.
\bibitem{Joyce9} D.~D.~Joyce, {\it Special Lagrangian Submanifolds
with Isolated Conical Singularities. V.~Survey and Applications}, J.~Differential Geom. {\bf 63}, 299-347, 2003.
\bibitem{Joyce1} D. D. Joyce, {\it Riemannian Holonomy Groups and Calibrated Geometry}, 
Oxford Graduate Texts in Mathematics {\bf 12}, OUP, Oxford, 2007. 
\bibitem{JoySalur} D.~D.~Joyce and S.~Salur, {\it Deformations of Asymptotically Cylindrical Coassociative
Submanifolds with Fixed Boundary}, Geom.~Topol. {\bf 9}, 1115-1146, 2005 (electronic).
\bibitem{LotayKov} A.~Kovalev and J.~D.~Lotay, {\it Deformations of Compact Coassociative 4-folds with Boundary}, 
 J.~Geom.~Phys. {\bf 59}, 63-73, 2009.
\bibitem{Lang} S.~Lang, {\it Differential Manifolds},
Addison-Wesley, Reading, Massachusetts, 1972.
\bibitem{Lang2} S.~Lang, {\it Real Analysis}, Addison-Wesley, Reading, Massachusetts, 1983.
\bibitem{Lockhart} R.~B.~Lockhart, {\it Fredholm, Hodge and Liouville
Theorems on Noncompact Manifolds}, Trans.~Amer.~Math.~Soc. {\bf 301}, 1-35, 1987.
\bibitem{LockhartMcOwen} R.~B.~Lockhart and R.~C.~McOwen, {\it
Elliptic Differential Operators on Noncompact
Manifolds},  Ann.~Sc.~Norm.~Super.~Pisa Cl.~Sci. {\bf 12}, 409-447, 1985.
\bibitem{Lotay2} J.~Lotay, {\it 2-Ruled Calibrated 4-folds in
$\R^7$ and $\R^8$}, J.~London Math.~Soc. {\bf 74}, 219-243, 2006.
\bibitem{Lotaysing} J.~D.~Lotay, {\it Coassociative 4-folds with Conical Singularities}, Comm. Anal. Geom. {\bf 15}, 891-946, 2008.
\bibitem{Lotaydesing} J.~D.~Lotay, {\it Desingularization of Coassociative 4-folds with
Conical Singularities}, Geom.~Funct.~Anal. {\bf 18}, 2055--2100, 2008.
\bibitem{Lotaylag} J.~D.~Lotay, {\it Ruled Lagrangian Submanifolds of the 6-Sphere}, to appear in Trans.~Amer.~Math.~Soc.
\bibitem{Marshall} S.~P.~Marshall, {\it Deformations of Special
Lagrangian Submanifolds}, DPhil thesis, Oxford University, Oxford, 2002.
\bibitem{Mazya} V.~G.~Maz'ya and B.~Plamenevskij, {\it Elliptic Boundary Value Problems}, Amer.~Math.~Soc.~Transl. {\bf 123}, 
1--56, 1984.
\bibitem{McLean} R.~C.~McLean, {\it Deformations of Calibrated
Submanifolds}, Comm. Anal. Geom. {\bf 6},
705-747, 1998.
\bibitem{Pacini} T.~Pacini, {\it Deformations of Asymptotically
Conical Special Lagrangian Submanifolds}, Pacific J.~Math. {\bf
215}, 151-181, 2004.
\bibitem{Salur} S.~Salur, {\it Deformations of Asymptotically Cylindrical Coassociative Submanifolds with Moving Boundary},
preprint, arXiv:math.DG/0601420.
\end{thebibliography}
\end{document}